\documentclass[final,3p]{elsarticle}
\usepackage{amsfonts}
\usepackage{graphicx}
\usepackage{amsmath}
\usepackage{amssymb}
\usepackage{amsthm}
\usepackage{mathrsfs}
\usepackage{color}
\usepackage{makecell,multirow,diagbox}
\usepackage{amsmath}
\usepackage{amsfonts}
\usepackage{amsthm}
\usepackage{amsmath,amsxtra,latexsym,amsthm, amssymb, amscd, pb-diagram}
\usepackage{subfigure}

\newtheorem{dn}{Definition}[section]
\newtheorem{bdt}{Inequality}[section]
\newtheorem{dl}{Theorem}[section]
\newtheorem{md}{Proposition}[section]
\newtheorem{bd}{Lemma}[section]
\newtheorem{hq}{Corollary}[section]
\newtheorem{nx}{Remark}[section]
\newtheorem{vd}{Example}[section]
\newtheorem{tc}{Property}[section]
\newtheorem{cm}{Proof}[section]
\newcommand{\R}{\mathbb{R}}

\newcommand{\Z}{\mathbb{Z}}

\newcommand{\e}{\varepsilon}
\newcommand{\ity}{\infty}
\newcommand{\f}{\frac}

\newcommand{\bbd}{\begin{bd}}
\newcommand{\ebd}{\end{bd}}
\newcommand{\bbdt}{\begin{bdt}}
\newcommand{\ebdt}{\end{bdt}}
\newcommand{\bdn}{\begin{dn}}
\newcommand{\edn}{\end{dn}}
\newcommand{\bhq}{\begin{hq}}
\newcommand{\ehq}{\end{hq}}
\newcommand{\bdl}{\begin{dl}}
\newcommand{\edl}{\end{dl}}
\newcommand{\bnx}{\begin{nx}}
\newcommand{\enx}{\end{nx}}
\newcommand{\bmd}{\begin{md}}
\newcommand{\emd}{\end{md}}
\newcommand{\bcm}{\begin{cm}}
\newcommand{\ecm}{\end{cm}}
\newcommand{\bvd}{\begin{vd}}
\newcommand{\evd}{\end{vd}}
\newcommand{\btc}{\begin{tc}}
\newcommand{\etc}{\end{tc}}

\usepackage{epsfig}
\usepackage{enumitem}
\journal{Journal of Mathematical Analysis and Applications}

\begin{document}

\begin{frontmatter}




\author{Tuan Anh Dao\fnref{label1,label2}}
\ead{daotuananh.fami@gmail.com}
\author{Michael Reissig\corref{cor2}\fnref{label2}}
\ead{reissig@math.tu-freiberg.de}

%\fntext[a]{School of Applied Mathematics and Informatics, Hanoi University of Science and Technology, No.1 Dai Co Viet road, Hanoi, Vietnam}
%\cortext[cor1]{}
%\fntext[b]{Faculty for Mathematics and Computer Science, TU Bergakademie Freiberg, Pr\"{u}ferstr. 9, 09596, Freiberg, Germany}

\cortext[cor2]{Corresponding author.}

% \address{Address\fnref{label3}}

\title{An application of $L^1$ estimates for oscillating integrals to parabolic like semi-linear structurally damped $\sigma$-evolution models}

% \author[label1,label2]{Tuan Anh Dao}

\address[label1]{School of Applied Mathematics and Informatics, Hanoi University of Science and Technology, No.1 Dai Co Viet road, Hanoi, Vietnam}
\address[label2]{Faculty for Mathematics and Computer Science, TU Bergakademie Freiberg, Pr\"{u}ferstr. 9, 09596, Freiberg, Germany}

% \author{}
% \address{}

\begin{abstract}
We study the following Cauchy problems for semi-linear structurally damped $\sigma$-evolution models:
\begin{equation*}
u_{tt}+ (-\Delta)^\sigma u+ \mu (-\Delta)^\delta u_t = f(u,u_t),\, u(0,x)= u_0(x),\, u_t(0,x)=u_1(x)
\end{equation*}
with $\sigma \ge 1$, $\mu>0$ and $\delta \in (0,\frac{\sigma}{2})$. Here the function $f(u,u_t)$ stands for the power nonlinearities $|u|^{p}$ and $|u_t|^{p}$ with a given number $p>1$. We are interested in investigating $L^{1}$ estimates for oscillating integrals in the presentation of the solutions to the corresponding linear models with vanishing right-hand sides by applying the theory of modified Bessel functions and Fa\`{a} di Bruno's formula. By assuming additional $L^{m}$ regularity on the initial data, we use $(L^{m}\cap L^{q})- L^{q}$ and $L^{q}- L^{q}$ estimates with $q\in (1,\ity)$ and $m\in [1,q)$, to prove the global (in time) existence of small data Sobolev solutions to the above semi-linear models from suitable function spaces basing on $L^q$ spaces.
\end{abstract}


\end{frontmatter}


%=================================================================================={Introduction}
\section{Introduction}
\label{Introduction}
In this paper, we consider the following two Cauchy problems:
\begin{equation}
 u_{tt}+ (-\Delta)^\sigma u+ \mu (-\Delta)^{\delta} u_t=|u|^p,\, u(0,x)= u_0(x),\, u_t(0,x)=u_1(x) \label{pt6.1}
\end{equation}
and
\begin{equation}
 u_{tt}+ (-\Delta)^\sigma u+ \mu (-\Delta)^{\delta} u_t=|u_t|^p,\, u(0,x)= u_0(x),\, u_t(0,x)=u_1(x) \label{pt6.2}
\end{equation}
with $\sigma \ge 1$, $\mu>0$ and $\delta \in (0,\frac{\sigma}{2})$. The corresponding linear models with vanishing right-hand side are
\begin{equation}
 u_{tt}+ (-\Delta)^\sigma u+ \mu (-\Delta)^{\delta} u_t=0,\text{ } u(0,x)= u_0(x),\text{ } u_t(0,x)=u_1(x). \label{pt6.3}
\end{equation}
Recently, there are several papers (see, for instance, \cite{DabbiccoEbert2014,DabbiccoReissig,NarazakiReissig}) concerning the special case  $\sigma=1$ to the linear Cauchy problems (\ref{pt6.3}), that is, to
\begin{equation}
 u_{tt}-\Delta u+ \mu (-\Delta)^{\delta} u_t=0,\text{ } u(0,x)= u_0(x),\text{ } u_t(0,x)=u_1(x) \label{pt6.4}
\end{equation}
with $\delta \in (0,1]$. In particular, in \cite{NarazakiReissig} the authors divided the phase space into two parts including sufficiently small and sufficiently large frequencies in order to study Fourier multipliers with oscillations in the representation of the solution to (\ref{pt6.4}). More in detail, to do this there appeared two main strategies in \cite{NarazakiReissig}. They applied heavily radial symmetry combined with the theory of modified Bessel functions (see also \cite{ReissigEbert}) and took into considerations the connection to Fourier multipliers appearing for wave models, respectively, for small frequencies and large frequencies. Consequently, having $L^1$ estimates for oscillating integrals was to conclude $L^p- L^q$ estimates not necessarily on the conjugate line for the solutions to (\ref{pt6.4}). In the case of semi-linear structurally damped wave models (\ref{pt6.1}) with $\sigma=1$ and $\delta \in (0,1]$ (see \cite{DabbiccoReissig}), the authors proved the global (in time) existence of small data solutions in low space dimensions by using classical energy estimates, i.e., estimates on the base of $L^2$ norms. In addition, in \cite{DabbiccoEbert2014} some suitable high frequencies $L^q- L^q$ estimates for the solution to (\ref{pt6.4}) have been obtained for $\delta \in (0,\frac{1}{4})$. Meanwhile, in the remaining case $\sigma \in [\frac{1}{4},1)$ the authors developed these estimates relying on some techniques in \cite{NarazakiReissig}. Then, some global (in time) existence results of small data solutions were presented for ``parabolic like models'' corresponding to (\ref{pt6.1}) with $\sigma=1$ and $\delta \in (0,\frac{1}{2})$.\medskip

More recently, the use of $(L^{1}\cap L^{2})- L^{2}$ estimates to (\ref{pt6.3}), i.e., the mixing of additional $L^{1}$ regularity for the data on the basis of $L^{2}- L^{2}$ estimates was investigated in \cite{DuongKainaneReissig} to study semi-linear $\sigma$-evolution models (\ref{pt6.1}) and (\ref{pt6.2}) with $\delta= \frac{\sigma}{2}$. The effective tools that the authors applied were results from Harmonic Analysis such as Gagliardo-Nirenberg inequality, fractional powers and embeddings into $L^{\ity}$ (see also \cite{Palmierithesis}). Some classical versions of the Gagliardo-Nirenberg inequality can be found, for example, in \cite{DabbiccoReissig,Ozawa,Bui}. Moreover, another approach in \cite{DabbiccoEbert} was to derive sharp $L^{p}- L^{q}$ estimates, with $1<p \le q<\ity$, to the linear models (\ref{pt6.3}) and some $L^{q}$ estimates for the solutions and some of their derivatives, with $q \in (1,\ity)$, to the semi-linear models (\ref{pt6.1}) and (\ref{pt6.2}) in the case with $\delta \in [0,\frac{\sigma}{2}]$. In particular, here the authors found an explicit way to obtain these estimates by using the Mikhlin-H\"{o}rmander multiplier theorem for kernels localized at high frequencies. Due to the lack of $L^{1}- L^{1}$ estimates, they used two different strategies to look for the global (in time) existence of small data solutions to semi-linear models. On the one hand, they took account of additional $L^{1}\cap L^{\ity}$ regularity in the first case with $\delta= \frac{\sigma}{2}$. Additional $L^{\eta}\cap L^{\bar{q}}$ regularity, on the other hand, was replaced for any small $\eta$ and large $\bar{q}$ in the second case with $\delta \in (0,\frac{\sigma}{2})$.\medskip

The motivation of this paper is to derive $L^{p}- L^{q}$ estimates for solutions to (\ref{pt6.3}) with $1 \le p \le q \le \ity$. Then, by mixing additional $L^{m}$ regularity for the data on the basis of $L^{q}- L^{q}$ estimates with $1 \le m< q< \ity$, we prove the global (in time) existence of small data solutions to semi-linear models (\ref{pt6.1}) and (\ref{pt6.2}) as well. For this reason, the first main goal of the present paper is to get $L^{1}$ estimates for oscillating integrals in the presentation of solutions to (\ref{pt6.3}) by using the theory of modified Bessel functions and Fa\`{a} di Bruno's formula (see, for instance, \cite{FrancescoBruno,Bui}). It is reasonable to apply Fa\`{a} di Bruno's formula since the connection to Fourier multipliers appearing for wave models fails to $\sigma$-evolution models. The second main goal of this paper is to use different strategies allowing no loss of decay and some loss of decay combined with loss of regularity to deal with semi-linear problems (\ref{pt6.1}) and (\ref{pt6.2}). \medskip

{\it Loss of regularity} (see, for example, \cite{CicognaniHirosawaReissig,DabbiccoEbert,Miyachi1980,Peral}) is a well-known phenomenon describing the effect that the regularity of the obtained solutions to semi-linear models is less than those of the initial data. This phenomenon appearing in our global (in time) existence results is due to the singular behavior of time-dependent coefficients in the estimates of solutions to the linear models localized to high frequencies as $t\longrightarrow +0$. However, we can compensate this difficulty by assuming higher regularity for the data. \\

{\it Loss of decay} is understood when the decay rates in estimates of the solutions to semi-linear models are worse than those given for the solutions to the linear models with vanishing right-hand side. Additional benefits of allowing loss of decay (see \cite{DabbiccoEbert2014}) are to show how the restrictions to the admissible exponents $p$ could be relaxed. In this paper, some new tools from Harmonic Analysis (see \cite{Palmierithesis}) play an important role to prove our global (in time) existence results.\medskip

The scheme of this paper is organized as follows:
\begin{itemize}[leftmargin=*]
\item In Section $2$, we state the main results. In particular, in Section $2.1$ we state $L^{p}- L^{q}$ estimates not necessarily on the conjugate line for solutions to (\ref{pt6.3}). In Section $2.2$, we state our global (in time) existence results of small data solutions to (\ref{pt6.1}) and (\ref{pt6.2}) without loss of decay and with loss of decay combined with loss of regularity.
\item In Section $3$, we present estimates for the solutions to (\ref{pt6.3}). We devote to the proof of $L^{1}$ estimates, $L^{\ity}$ estimates and $L^{r}$ estimates, respectively, in Sections $3.1$, $3.2$ and $3.3$. Finally, in Section $3.4$ we state $L^{q}- L^{q}$ estimates assuming additional $L^{m}$ regularity for the data with $q\in (1,\ity)$ and $m\in [1,q)$.
\item In Section $4$, we prove our global (in time) existence results to (\ref{pt6.1}) and (\ref{pt6.2}).
\item In Section $5$, we state some concluding remarks and open problems.
\end{itemize}

For the ease of reading, in this paper we use the following notations.\medskip

\noindent\textbf{Notation 1.} We write $f\lesssim g$ when there exists a constant $C>0$ such that $f\le Cg$, and $f \approx g$ when $g\lesssim f\lesssim g$.\medskip

\noindent\textbf{Notation 2.} $H^{a,q}$ and $\dot{H}^{a,q}$, with $a \ge 0$ and $q \in (1,\ity)$, denote Bessel and Riesz potential spaces based on $L^q$. The abbreviations $\big<D\big>^{a}$ and $|D|^{a}$ stand for the pseudo-differential operators with symbols $\big<\xi\big>^{a}$ and $|\xi|^{a}$, respectively.\medskip

\noindent\textbf{Notation 3.} We denote $[s]^+:= \max\{s,0\}$ as the positive part of $s \in \R$, and $\lceil s \rceil:= \min \big\{k \in \Z \,\, : \,\, k\ge s \big\}$.

\noindent\textbf{Notation 4.} We introduce the spaces $\mathcal{A}^{s}_{m,q}:= \big(L^m \cap H^{s,q}\big) \times \big(L^m \cap H^{[s-\sigma]^+,q}\big)$ with the norm
$$\|(u_0,u_1)\|_{\mathcal{A}^{s}_{m,q}}:=\|u_0\|_{L^m}+ \|u_0\|_{H^{s,q}}+ \|u_1\|_{L^m}+ \|u_1\|_{H^{[s-\sigma]^+,q}}. $$
Here $s\ge 0$, $q \in (1,\infty)$ and $m \in [1,q)$.

\noindent\textbf{Notation 5.} We fix the constants $s_0:= \big(2+\big[\frac{n}{2}\big]\big)(\sigma-2\delta)$, $n_0:= \frac{6\delta-2\sigma}{\sigma-2\delta}$ and $n_1:=\f{4mq(\sigma-\delta)}{q-m}$.

%=================================================================================={Main results}
\section{Main results}
\label{Main results}
\subsection{$L^p- L^q$ estimates not necessarily on the conjugate line}
\noindent\textbf{Theorem 1.}\textit{ Let $\delta\in \big(0,\f{\sigma}{2}\big)$ in (\ref{pt6.3}) and $1\le p\le q\le \ity$. Then, the solutions to (\ref{pt6.3}) satisfy the following $L^p- L^q$ estimates:
\begin{align*}
\big\||D|^a u(t,\cdot)\big\|_{L^q}& \lesssim
\begin{cases}
t^{-(2+[\frac{n}{2}])(\frac{\sigma}{2\delta}-1)\frac{1}{r} -\frac{n}{2\delta}(1-\frac{1}{r})-\frac{a}{2\delta}}\|u_0\|_{L^p}+ t^{1-(1+[\frac{n}{2}])(\frac{\sigma}{2\delta}-1)\frac{1}{r} -\frac{n}{2\delta}(1-\frac{1}{r})-\frac{a}{2\delta}}\|u_1\|_{L^p} \text{ if } t\in (0,1], &\\
t^{-\frac{n}{2(\sigma-\delta)}(1-\frac{1}{r})-\frac{a}{2(\sigma-\delta)}} \|u_0\|_{L^p}+ t^{1-\frac{n}{2(\sigma-\delta)}(1-\frac{1}{r})-\frac{a}{2(\sigma-\delta)}}\|u_1\|_{L^p} \text{ if } t\in[1,\ity), &
\end{cases}\\
\big\||D|^a u_t(t,\cdot)\big\|_{L^q}& \lesssim
\begin{cases}
t^{1-(1+[\frac{n}{2}])(\frac{\sigma}{2\delta}-1)\frac{1}{r} -\frac{n}{2\delta}(1-\frac{1}{r})-\frac{a+2\sigma}{2\delta}}\|u_0\|_{L^p}+ t^{-(2+[\frac{n}{2}])(\frac{\sigma}{2\delta}-1)\frac{1}{r} -\frac{n}{2\delta}(1-\frac{1}{r})-\frac{a}{2\delta}}\|u_1\|_{L^p} \text{ if } t\in (0,1], &\\
t^{-\frac{n}{2(\sigma-\delta)}(1-\frac{1}{r})-\frac{a+2\delta}{2(\sigma-\delta)}} \|u_0\|_{L^p}+ t^{1-\frac{n}{2(\sigma-\delta)}(1-\frac{1}{r})-\frac{a+2\delta}{2(\sigma-\delta)}}\|u_1\|_{L^p} \text{ if } t\in[1,\ity), &
\end{cases}
\end{align*}
where $1+ \frac{1}{q}= \frac{1}{r}+\frac{1}{p}$, for any non-negative number $a$ and for all $n \ge 1$.}

\subsection{Global (in time) existence of small data solutions}
In the following statements we use $s_0:= \big(2+\big[\frac{n}{2}\big]\big)(\sigma-2\delta)$, $n_0:= \frac{6\delta-2\sigma}{\sigma-2\delta}$ and $n_1:=\f{4mq(\sigma-\delta)}{q-m}$.\medskip

\noindent In the first case, we obtain solutions to (\ref{pt6.1}) from energy space on the base of $L^q$.\medskip

\noindent\textbf{Theorem 2-A.}\textit{ Let $q \in (1,\ity)$ be a fixed constant and $m\in [1,q)$. We assume the conditions $\big[\frac{n}{2}\big]< n_0$ and
\begin{equation}
p> 1+\frac{\max\big\{n-\frac{m}{q}n+m\sigma,\, 4m(\sigma-\delta)\big\}}{n-2m(\sigma-\delta)}. \label{exponent2A}
\end{equation}
Moreover, we suppose the following conditions:
\begin{equation}
p\in \Big[\frac{q}{m},\ity \Big)  \text{ if }  n\le q\sigma, \text{ or } p \in \Big[\f{q}{m}, \f{n}{n-q\sigma}\Big]  \text{ if } n \in \Big(q\sigma, \f{q^2 \sigma}{q-m}\Big]. \label{GN2A}
\end{equation}
Then, there exists a constant $\e>0$ such that for any small data
$$(u_0,u_1) \in \mathcal{A}^{\sigma+s_0}_{m,q} \text{ satisfying the assumption } \|(u_0,u_1)\|_{\mathcal{A}^{\sigma+s_0}_{m,q}} \le \e, $$
we have a uniquely determined global (in time) small data energy solution (on the base of $L^q$) \[ u\in C([0,\ity),H^{\sigma,q})\cap C^1([0,\ity),L^q)\]
to (\ref{pt6.1}). The following estimates hold:
\begin{align}
\|u(t,\cdot)\|_{L^q}& \lesssim (1+t)^{1-\frac{n}{2(\sigma-\delta)}(1-\frac{1}{r})} \|(u_0,u_1)\|_{\mathcal{A}^{\sigma+s_0}_{m,q}}, \label{decayrate2A1}\\
\big\||D|^\sigma u(t,\cdot)\big\|_{L^q}& \lesssim (1+t)^{1-\frac{n}{2(\sigma-\delta)}(1-\frac{1}{r})-\frac{\sigma}{2(\sigma-\delta)}} \|(u_0,u_1)\|_{\mathcal{A}^{\sigma+s_0}_{m,q}}, \label{decayrate2A2}\\
\|u_t(t,\cdot)\|_{L^q}& \lesssim (1+t)^{1-\frac{n}{2(\sigma-\delta)}(1-\frac{1}{r})-\frac{\delta}{\sigma-\delta}} \|(u_0,u_1)\|_{\mathcal{A}^{\sigma+s_0}_{m,q}}, \label{decayrate2A3}
\end{align}
where $1+\frac{1}{q}=\frac{1}{r}+\frac{1}{m}$.}\medskip

\noindent\textbf{Theorem 2-B.}\textit{ Under the assumptions of Theorem 2-A, if condition (\ref{exponent2A}) is replaced by $n> n_1$, then we have the same conclusions of Theorem 2-A. But the estimates (\ref{decayrate2A1})-(\ref{decayrate2A3}) are modified in the following way:
\begin{align}
\big\|\big(|D|^\sigma u(t,\cdot), u(t,\cdot)\big)\big\|_{L^q}& \lesssim (1+t)^{\frac{1}{p}\big(1-\frac{n}{2(\sigma-\delta)}(1-\frac{1}{r})\big)} \|(u_0,u_1)\|_{\mathcal{A}^{\sigma+s_0}_{m,q}}, \label{decayrate2B1}\\
\|u_t(t,\cdot)\|_{L^q}& \lesssim (1+t)^{1-\frac{n}{2(\sigma-\delta)}(1-\frac{1}{r})} \|(u_0,u_1)\|_{\mathcal{A}^{\sigma+s_0}_{m,q}}. \label{decayrate2B2}
\end{align}}
In the second case, we obtain Sobolev solutions to (\ref{pt6.1}).\medskip

\noindent\textbf{Theorem 3-A.}\textit{ Let $q \in (1,\ity)$ be a fixed constant, $m\in [1,q)$ and $0 < s < \sigma$. We assume the conditions $\big[\frac{n}{2}\big]< n_0$ and
\begin{equation}
p> 1+\frac{\max\big\{n-\frac{m}{q}n+ms,\, 4m(\sigma-\delta)\big\}}{n-2m(\sigma-\delta)}. \label{exponent3A}
\end{equation}
Moreover, we suppose the following conditions:
\begin{equation}
p\in \Big[\frac{q}{m},\ity \Big)  \text{ if } n\le qs, \text{ or }p \in \Big[\f{q}{m}, \f{n}{n-qs}\Big] \text{ if } n \in \Big(qs, \f{q^2 s}{q-m}\Big]. \label{GN3A}
\end{equation}
Then, there exists a constant $\e>0$ such that for any small data
$$(u_0,u_1) \in \mathcal{A}^{s+s_0}_{m,q} \text{ satisfying the assumption } \|(u_0,u_1)\|_{\mathcal{A}^{s+s_0}_{m,q}}\le \e, $$
we have a uniquely determined global (in time) small data Sobolev solution \[ u\in C([0,\ity),H^{s,q})\]
to (\ref{pt6.1}). The following estimates hold:
\begin{align}
\|u(t,\cdot)\|_{L^q}& \lesssim (1+t)^{1-\frac{n}{2(\sigma-\delta)}(1-\frac{1}{r})} \|(u_0,u_1)\|_{\mathcal{A}^{s+s_0}_{m,q}}, \label{decayrate3A1}\\
\big\||D|^s u(t,\cdot)\big\|_{L^q}& \lesssim (1+t)^{1-\frac{n}{2(\sigma-\delta)}(1-\frac{1}{r})-\frac{s}{2(\sigma-\delta)}} \|(u_0,u_1)\|_{\mathcal{A}^{s+s_0}_{m,q}}, \label{decayrate3A2}
\end{align}
where $1+\frac{1}{q}=\frac{1}{r}+\frac{1}{m}$.}\medskip

\noindent\textbf{Theorem 3-B.}\textit{ Under the assumptions of Theorem 3-A, if the condition (\ref{exponent3A}) is replaced by $n> n_1$, then we have the same conclusions of Theorem 3-A. But the estimates (\ref{decayrate3A1})-(\ref{decayrate3A2}) are modified in the following way:
\begin{equation}
\big\|\big(|D|^s u(t,\cdot), u(t,\cdot)\big)\big\|_{L^q} \lesssim (1+t)^{\frac{1}{p}\big(1-\frac{n}{2(\sigma-\delta)}(1-\frac{1}{r})\big)} \|(u_0,u_1)\|_{\mathcal{A}^{s+s_0}_{m,q}}. \label{decayrate3B}
\end{equation}}
In the third case, we obtain solutions to (\ref{pt6.1}) belonging to the energy space (on the base of $L^q$) with a suitable higher regularity.\medskip

\noindent\textbf{Theorem 4-A.}\textit{ Let $q \in (1,\ity)$ be a fixed constant, $m\in [1,q)$ and $\sigma<s\le \sigma+\frac{n}{q}$. We assume that the exponent $p$ satisfies the conditions $p>1+ \lceil s-\sigma \rceil$ and
\begin{equation}
p> 1+\frac{\max\big\{n-\frac{m}{q}n+ms,\, 4m(\sigma-\delta)\big\}}{n-2m(\sigma-\delta)}, \label{exponent4A}
\end{equation}
where $\big[\frac{n}{2}\big]< n_0$. Moreover, we suppose the following conditions:
\begin{equation}
p \in \Big[\f{q}{m}, \ity \Big) \text{ if } n\le qs, \text{ or } p \in \Big[\f{q}{m}, 1+\f{q\sigma}{n-qs}\Big] \text{ if } n \in \Big(qs, qs+\f{qm\sigma}{q-m}\Big]. \label{GN4A}
\end{equation}
Then, there exists a constant $\e>0$ such that for any small data
$$(u_0,u_1) \in \mathcal{A}^{s+s_0}_{m,q}\text{ satisfying the assumption }\|(u_0,u_1)\|_{\mathcal{A}^{s+s_0}_{m,q}}\le \e, $$
we have a uniquely determined global (in time) small data energy solution \[ u\in C([0,\ity),H^{s,q})\cap C^1([0,\ity),H^{s-\sigma,q})\]
to (\ref{pt6.1}). The following estimates hold:
\begin{align}
\|u(t,\cdot)\|_{L^q}& \lesssim (1+t)^{1-\frac{n}{2(\sigma-\delta)}(1-\frac{1}{r})} \|(u_0,u_1)\|_{\mathcal{A}^{s+s_0}_{m,q}}, \label{decayrate4A1}\\
\big\||D|^s u(t,\cdot)\big\|_{L^q}& \lesssim (1+t)^{1-\frac{n}{2(\sigma-\delta)}(1-\frac{1}{r})-\frac{s}{2(\sigma-\delta)}} \|(u_0,u_1)\|_{\mathcal{A}^{s+s_0}_{m,q}}, \label{decayrate4A2}\\
\|u_t(t,\cdot)\|_{L^q}& \lesssim (1+t)^{1-\frac{n}{2(\sigma-\delta)}(1-\frac{1}{r})-\frac{\delta}{\sigma-\delta}} \|(u_0,u_1)\|_{\mathcal{A}^{s+s_0}_{m,q}}, \label{decayrate4A3}\\
\big\||D|^{s-\sigma} u_t(t,\cdot)\big\|_{L^q}& \lesssim (1+t)^{1-\frac{n}{2(\sigma-\delta)}(1-\frac{1}{r})-\frac{s-\sigma+2\delta}{2(\sigma-\delta)}} \|(u_0,u_1)\|_{\mathcal{A}^{s+s_0}_{m,q}}, \label{decayrate4A4}
\end{align}
where $1+\frac{1}{q}=\frac{1}{r}+\frac{1}{m}$.}\medskip

\noindent\textbf{Theorem 4-B.}\textit{ Under the assumptions of Theorem 4-A, if condition (\ref{exponent4A}) is replaced by $n> n_1$, then we have the same conclusions of Theorem 4-A. But the estimates (\ref{decayrate4A1})-(\ref{decayrate4A4}) are modified in the following way:
\begin{align}
\big\|\big(|D|^s u(t,\cdot), u(t,\cdot)\big)\big\|_{L^q}& \lesssim (1+t)^{\frac{1}{p}\big(1-\frac{n}{2(\sigma-\delta)}(1-\frac{1}{r})\big)} \|(u_0,u_1)\|_{\mathcal{A}^{s+s_0}_{m,q}}, \label{decayrate4B1}\\
\big\|\big(|D|^{s-\sigma} u_t(t,\cdot), u_t(t,\cdot)\big)\big\|_{L^q}& \lesssim (1+t)^{1-\frac{n}{2(\sigma-\delta)}(1-\frac{1}{r})} \|(u_0,u_1)\|_{\mathcal{A}^{s+s_0}_{m,q}}. \label{decayrate4B2}
\end{align}}
Finally, we obtain large regular solutions to (\ref{pt6.1}) by using the fractional Sobolev embedding.\medskip

\noindent\textbf{Theorem 5-A.}\textit{ Let $s >\sigma+\frac{n}{q}$. Let $q \in (1,\ity)$ be a fixed constant and $m\in [1,q)$. We assume that the exponent $p$ satisfies the conditions $p> 1+ s-\sigma$ and
\begin{equation}
p> 1+\frac{\max\big\{n-\frac{m}{q}n+ms,\, 4m(\sigma-\delta)\big\}}{n-2m(\sigma-\delta)}, \label{exponent5A}
\end{equation}
where $\big[\frac{n}{2}\big]< n_0$. Moreover, we suppose the following conditions:
\begin{equation}
p \in \Big[\f{q}{m},\ity\Big) \text{ and } n> 2m(\sigma-\delta). \label{GN5A}
\end{equation}
Then, there exists a constant $\e>0$ such that for any small data
$$(u_0,u_1) \in \mathcal{A}^{s+s_0}_{m,q}\text{ satisfying the assumption }\|(u_0,u_1)\|_{\mathcal{A}^{s+s_0}_{m,q}}\le \e, $$
we have a uniquely determined global (in time) small data energy solution \[ u\in C([0,\ity),H^{s,q})\cap C^1([0,\ity),H^{s-\sigma,q})\]
to (\ref{pt6.1}). Moreover, the estimates (\ref{decayrate4A1})-(\ref{decayrate4A4}) hold.}\medskip

\noindent\textbf{Theorem 5-B.}\textit{ Under the assumptions of Theorem 5-A, if the condition (\ref{exponent5A}) is replaced by $n> n_1$, then we have same conclusions of Theorem 5-A. But the estimates (\ref{decayrate4A1})-(\ref{decayrate4A4}) are modified in the following way:
\begin{align}
\big\|\big(|D|^s u(t,\cdot), u(t,\cdot)\big)\big\|_{L^q}& \lesssim (1+t)^{\frac{1}{p}\big(1-\frac{n}{2(\sigma-\delta)}(1-\frac{1}{r})\big)} \|(u_0,u_1)\|_{\mathcal{A}^{s+s_0}_{m,q}}, \label{decayrate5B1}\\
\big\|\big(|D|^{s-\sigma} u_t(t,\cdot), u_t(t,\cdot)\big)\big\|_{L^q}& \lesssim (1+t)^{1-\frac{n}{2(\sigma-\delta)}(1-\frac{1}{r})} \|(u_0,u_1)\|_{\mathcal{A}^{s+s_0}_{m,q}}. \label{decayrate5B2}
\end{align}}
Finally, we obtain large regular solutions to (\ref{pt6.2}) by using the fractional Sobolev embedding. \medskip

\noindent\textbf{Theorem 6-A.}\textit{ Let $s >\sigma+\frac{n}{q}$. Let $q \in (1,\ity)$ be a fixed constant and $m\in [1,q)$. We assume that the exponent $p$ satisfies the conditions $p> 1+ s-\sigma$ and
\begin{equation}
p> 1+\frac{\max\big\{n-\frac{m}{q}n+m(s-2\delta),\, 2m(2\sigma-3\delta)\big\}}{n-2m(\sigma-2\delta)}, \label{exponent6A}
\end{equation}
where $\big[\frac{n}{2}\big]< n_0$. Moreover, we suppose the following conditions:
\begin{equation}
p \in \Big[\f{q}{m},\ity\Big) \text{ and } n> 2m(\sigma-2\delta). \label{GN6A}
\end{equation}
Then, there exists a constant $\e>0$ such that for any small data
$$(u_0,u_1) \in \mathcal{A}^{s+s_0}_{m,q}\text{ satisfying the assumption }\|(u_0,u_1)\|_{\mathcal{A}^{s+s_0}_{m,q}}\le \e, $$
we have a uniquely determined global (in time) small data energy solution \[ u\in C([0,\ity),H^{s,q})\cap C^1([0,\ity),H^{s-\sigma,q})\]
to (\ref{pt6.2}). Moreover, the estimates (\ref{decayrate4A1})-(\ref{decayrate4A4}) hold.}\medskip

\noindent\textbf{Theorem 6-B.}\textit{ Under the assumptions of Theorem 6-A, if the condition (\ref{exponent6A}) is replaced by $n> n_1$, then we have the same conclusions of Theorem 6-A. But the estimates (\ref{decayrate4A1})-(\ref{decayrate4A4}) are modified in the following way:
\begin{align}
\big\|\big(|D|^s u(t,\cdot), u(t,\cdot)\big)\big\|_{L^q}& \lesssim (1+t)^{1-\frac{n}{2(\sigma-\delta)}(1-\frac{1}{r})} \|(u_0,u_1)\|_{\mathcal{A}^{s+s_0}_{m,q}}, \label{decayrate6B1}\\
\big\|\big(|D|^{s-\sigma} u_t(t,\cdot), u_t(t,\cdot)\big)\big\|_{L^q}& \lesssim (1+t)^{\frac{1}{p}\big(1-\frac{n}{2(\sigma-\delta)}(1-\frac{1}{r})\big)} \|(u_0,u_1)\|_{\mathcal{A}^{s+s_0}_{m,q}}. \label{decayrate6B2}
\end{align}}

\begin{nx}
\fontshape{n}
\selectfont
There appears a loss of regularity $s_0$ of the solutions in all the above theorems with respect to the initial data. The phenomenon appears due to the application of some estimates for solutions to
 (\ref{pt6.3}) on the basis of $L^q$ with $q \in (1,\ity)$, to treat the semi-linear models (\ref{pt6.1}) and (\ref{pt6.2}). We will see this later in Proposition \ref{md3.7}.
\end{nx}

\begin{nx}
\fontshape{n}
\selectfont
For the estimates of the solutions to (\ref{pt6.1}) and (\ref{pt6.2}), in the above theorems A there appear the same decay rates as in the estimates for the solutions to (\ref{pt6.3}), i.e., no loss of decay appears. However, if we want to simplify some restrictions to the exponent $p$, for example, (\ref{exponent2A}), then we pay with further conditions for space dimension $n$, namely $n> n_1$ in the theorems B. Moreover, we can see that the decay rates for the solutions to the semi-linear models, for example in (\ref{decayrate2B1}) and (\ref{decayrate2B2}), are worse than those for solutions to the corresponding linear models, i.e., some loss of decay appears. This phenomenon is related to some of the used techniques in our proofs.
\end{nx}

\begin{nx}
\fontshape{n}
\selectfont
Let us compare our results with some known results. First, to the linear models (\ref{pt6.3}), in the special case of $\sigma=1$ one may show that the decay estimate for the solution itself appearing in Theorem 1 is almost the same as the corresponding one obtained in \cite{NarazakiReissig} if we consider the case of sufficiently large space dimensions $n$. If we set formally $\delta=\frac{\sigma}{2}$ in Theorem 1, then we can see that our results coincide with those in \cite{DuongKainaneReissig}. To the semi-linear models (\ref{pt6.1}), by putting $\sigma=1$, $q=2$ and $m=1$ we observe that the admissible exponents $p$ in Theorem 2-A are less flexible than those in \cite{DabbiccoReissig} for space dimensions $n=2,3,4$. However, in comparison with \cite{DabbiccoReissig} we want to point out that Theorem 2-A completely bring some flexibility for both $p$ and $n$, in general, due to the flexible choice of parameters $\sigma$, $\delta$, $q$ and $m$ (see also some of the examples below).
\end{nx}

\begin{vd}
\fontshape{n}
\selectfont
In the following examples, we choose $m=1$, $q=5$, $\sigma=2$ and $\delta=\frac{9}{10}$:
\begin{itemize}
\item If $n= 3$, then using Theorem 2-A we obtain $p \in \big(\frac{13}{2},\ity\big)$.
\item If $n= 3$ and $s= \frac{3}{2}$, then using Theorem 3-A we obtain $p \in \big(\frac{13}{2},\ity\big)$.
\item If $n= 3$ and $s= \frac{5}{2}$, then using Theorem 4-A we obtain $p \in \big(\frac{49}{8},\ity\big)$.
\item If $n= 5$ and $s= 5$, then using Theorem 5-A we obtain $p \in [5,\ity)$.
\item If $n= 3$ and $s= 5$, then using Theorem 6-A we obtain $p \in [5,\ity)$.
\end{itemize}
\end{vd}
\begin{vd}
\fontshape{n}
\selectfont
In the following examples, we choose $m=1$, $q=4$, $\sigma=2$ and $\delta=\frac{7}{8}$:
\begin{itemize}
\item If $n= 9$, then using Theorem 2-B we obtain $p \in [4,9]$.
\item If $n= 9$ and $s= \frac{9}{5}$, then using Theorem 3-B we obtain $p \in [4,5]$.
\item If $n= 9$ and $s= \frac{5}{2}$, then using Theorem 4-B we obtain $p \in [4,\ity)$.
\item If $n= 8$ and $s= 5$, then using Theorem 5-B we obtain $p \in \big(4,\ity)$.
\item If $n= 9$ and $s= 5$, then using Theorem 6-B we obtain $p \in \big(4,\ity)$.
\end{itemize}
\end{vd}

%=================================================================================={Linear estimates}
\section{Estimates for the solutions of the linear Cauchy problem}
\label{Linear estimates}
Using partial Fourier transformation to (\ref{pt6.3}), we obtain the following Cauchy problem for $v(t,\xi):=F_{x\rightarrow \xi}\big(u(t,x)\big)$, $v_0(\xi):=F_{x\rightarrow \xi}\big(u_0(x)\big)$ and $v_1(\xi):=F_{x\rightarrow \xi}\big(u_1(x)\big)$:
\begin{equation}
v_{tt}+ \mu |\xi|^{2\delta} v_t+ |\xi|^{2\sigma} v=0,\,\, v(0,\xi)= v_0(\xi),\,\, v_t(0,\xi)= v_1(\xi). \label{pt3.1}
\end{equation}
We can choose without loss of generality $\mu=1$ in (\ref{pt3.1}). The characteristic roots are
$$ \lambda_{1,2}=\lambda_{1,2}(\xi)= \f{1}{2}\Big(-|\xi|^{2\delta}\pm \sqrt{|\xi|^{4\delta}-4|\xi|^{2\sigma}}\Big). $$
The solution to (\ref{pt3.1}) is presented by the following formula (here we assume $\lambda_{1}\neq \lambda_{2}$):
\begin{equation}
v(t,\xi)= \frac{\lambda_1 e^{\lambda_2 t}-\lambda_2 e^{\lambda_1 t}}{\lambda_1- \lambda_2}v_0(\xi)+ \frac{e^{\lambda_1 t}-e^{\lambda_2 t}}{\lambda_1- \lambda_2}v_1(\xi)=: \hat{K_0}(t,\xi)v_0(\xi)+\hat{K_1}(t,\xi)v_1(\xi). \label{pt3.2}
\end{equation}
Taking account of the cases of small and large frequencies separately, we have
\begin{align}
&1.\,\, \lambda_1\sim -|\xi|^{2(\sigma- \delta)},\,\, \lambda_2\sim -|\xi|^{2\delta},\,\, \lambda_1-\lambda_2 \sim |\xi|^{2\delta} \text{ for small } |\xi|, \label{pt3.3}\\
&2.\,\, \lambda_{1,2} \sim -|\xi|^{2\delta}\pm i|\xi|^\sigma,\,\, \lambda_1-\lambda_2 \sim i|\xi|^\sigma \text{ for large } |\xi|. \label{pt3.4}
\end{align}
We now decompose the solution to (\ref{pt6.3}) into two parts localized separately to low and high frequencies, that is,
$$ u(t,x)= u_\chi(t,x)+ u_{1-\chi}(t,x), $$
where
$$u_\chi(t,x)=F^{-1}\big(\chi(|\xi|)v(t,\xi)\big) \text{ and } u_{1-\chi}(t,x)=F^{-1}\big(\big(1-\chi(|\xi|)\big)v(t,\xi)\big), $$
with a smooth cut-off function $\chi(|\xi|)$ equal to $1$ for small $|\xi|$ and vanishing for large $|\xi|$.

\subsection{$L^1$ estimates}
\subsubsection{Small frequencies}
Our approach is based on the paper \cite{NarazakiReissig}. According to the treatment of Propositions $4$ and $5$ in \cite{NarazakiReissig}, with minor modifications in the steps of the proofs we obtain the following $L^1$ estimates for small frequencies.
\bmd \label{md3.1.2}
The estimates
\begin{align*}
&\big\|F^{-1}\big(|\xi|^a \hat{K_0}\chi(|\xi|)\big)(t,\cdot)\big\|_{L^1}\lesssim
\begin{cases}
1 \text{ for } t\in (0,1], &\\
t^{-\frac{a}{2(\sigma-\delta)}} \text{ for } t\in[1,\ity),&
\end{cases}\\
&\big\|F^{-1}\big(|\xi|^a \hat{K_1}\chi(|\xi|)\big)(t,\cdot)\big\|_{L^1}\lesssim
\begin{cases}
t \text{ for } t\in (0,1], &\\
t^{1-\frac{a}{2(\sigma-\delta)}} \text{ for } t\in[1,\ity),&
\end{cases}
\end{align*}
hold for any non-negative number $a$.
\emd

\subsubsection{Large frequencies}
\bmd \label{md3.1.4}
The following estimates hold in $\R^n$ for any $n \ge 1$:
\begin{align*}
&\Big\| F^{-1}\Big(|\xi|^a \hat{K_0}\big(1-\chi(|\xi|)\big)\Big)(t,\cdot)\Big\|_{L^1}\lesssim \begin{cases}
t^{-(2+[\frac{n}{2}])(\frac{\sigma}{2\delta}-1)-\frac{a}{2\delta}} \text{ for } t\in (0,1], &\\
e^{-ct} \text{ for } t\in[1,\ity), &
\end{cases}\\
&\Big\| F^{-1}\Big(|\xi|^a \hat{K_1}\big(1-\chi(|\xi|)\big)\Big)(t,\cdot)\Big\|_{L^1}\lesssim \begin{cases}
t^{1-(1+[\frac{n}{2}])(\frac{\sigma}{2\delta}-1)-\frac{a}{2\delta}} \text{ for } t\in (0,1], &\\
e^{-ct} \text{ for } t\in[1,\ity), &
\end{cases}
\end{align*}
where $c$ is a suitable positive constant and $a$ is an arbitrary non-negative number $a$.
\emd
In order to obtain the desired estimates for the norm of the above Fourier multipliers localized to large frequencies, we can re-write
\begin{equation}
|\xi|^a \hat{K_0}(t,\xi) =e^{-\frac{1}{2}|\xi|^{2\delta}t}|\xi|^a \cos \Big( |\xi|^\sigma \sqrt{1-\frac{1}{4|\xi|^{2\sigma-4\delta}}} t \Big)+ e^{-\frac{1}{2} |\xi|^{2\delta}t}|\xi|^{a+2\delta}\frac{\sin \big( |\xi|^\sigma \sqrt{1-\frac{1}{4|\xi|^{2\sigma-4\delta}}}t \big)}{2|\xi|^\sigma \sqrt{1-\frac{1}{4|\xi|^{2\sigma-4\delta}}}}, \label{p3.2.1}
\end{equation}
and
\begin{equation}
|\xi|^a \hat{K_1}(t,\xi)=e^{-\frac{1}{2} |\xi|^{2\delta}t}|\xi|^a \frac{\sin \big(|\xi|^\sigma \sqrt{1-\frac{1}{4|\xi|^{2\sigma-4\delta}}}t \big)}{|\xi|^\sigma \sqrt{1-\frac{1}{4|\xi|^{2\sigma-4\delta}}}}. \label{p3.2.2}
\end{equation}
Hence, it seems to be reasonable to divide the proof into two steps. In the first step we derive $L^1$ estimates for the oscillating integrals
$$ F^{-1}\Big(e^{-c_1 |\xi|^{2\delta}t}|\xi|^{2\beta}\frac{\sin (c_2 |\xi|^\sigma t)}{|\xi|^\sigma}\big(1-\chi(|\xi|)\big)\Big), $$
and
$$ F^{-1}\Big(e^{-c_1 |\xi|^{2\delta}t}|\xi|^{2\beta}\cos (c_2 |\xi|^\sigma t) \big(1-\chi(|\xi|)\big)\Big), $$
where $\beta \ge 0$, $c_1$ is a positive constant and $c_2 \ne 0$ is a real constant. In the second step, we estimate the following two oscillating integrals:
$$ F^{-1}\Big(e^{-c_1 |\xi|^{2\delta}t}|\xi|^{2\beta}\frac{\sin \big(c_2 |\xi|^\sigma f(|\xi|) t\big)}{|\xi|^\sigma f(|\xi|)}\big(1-\chi(|\xi|)\big)\Big), $$
and
$$ F^{-1}\Big(e^{-c_1 |\xi|^{2\delta}t}|\xi|^{2\beta}\cos \big(c_2 |\xi|^\sigma f(|\xi|) t\big) \big(1-\chi(|\xi|)\big)\Big), $$
where $$f(|\xi|)=\sqrt{1-\frac{1}{4|\xi|^{2\sigma-4\delta}}}. $$

\bbd \label{bd3.1.3}
The following estimates hold in $\R^n$ for any $n \ge 1$:
$$\Big\|
F^{-1}\Big(e^{-c_1 |\xi|^{2\delta}t}|\xi|^{2\beta}\frac{\sin (c_2 |\xi|^\sigma t)}{|\xi|^\sigma}\big(1-\chi(|\xi|)\big)\Big)(t,\cdot)\Big\|_{L^1}
\lesssim \begin{cases}
t^{-(2+[\frac{n}{2}])(\frac{\sigma}{2\delta}-1)+\frac{\sigma-2\beta}{2\delta}} \text{ for } t\in (0,1], &\\
e^{-ct} \text{ for } t\in[1,\ity), &
\end{cases}$$
with $\beta \ge 0$ and $c$ is a suitable positive constant. Here $c_1$ is a positive and $c_2 \ne 0$ is a real constant.
\ebd
\begin{proof} We follow ideas from Proposition $4$ in \cite{NarazakiReissig}. Many steps in our proof are similar to those from Proposition $4$ devoting to small frequencies, nevertheless we will present the proof in detail to feel changes related to our interest for large frequencies. Let us divide the proof into two cases: $t\in [1,\ity)$ and $t\in (0,1]$. First, in order to treat the first case $t\in [1,\ity)$, we localize to small $|x|\le 1$. Then we obtain immediately the exponential decay. For this reason, we assume now $|x| \ge 1$. We introduce the function
$$I(t,x):=F^{-1}\Big(e^{-c_1 |\xi|^{2\delta}t}|\xi|^{2\beta-\sigma} \sin (c_2 |\xi|^\sigma t)\big(1-\chi(|\xi|)\big)\Big)(t,x). $$
Since the functions in the parenthesis are radially symmetric with respect to $\xi$, the inverse Fourier transform is radially symmetric with respect to $x$, too. Using modified Bessel functions we get
\begin{equation}
I(t,x)= c\int_0^\ity e^{-c_1 r^{2\delta}t}r^{2\beta-\sigma}\sin (c_2 r^\sigma t)\big(1-\chi(r)\big)r^{n-1} \tilde{J}_{\frac{n}{2}-1}(r|x|) dr. \label{l3.1.1}
\end{equation}
Let us consider odd spatial dimension $n=2m+1, m \ge 1.$ By introducing the vector field $Xf(r):= \frac{d}{dr}\big(\frac{1}{r}f(r)\big)$ as in \cite{NarazakiReissig} we carry out $m+1$ steps of partial integration to obtain
\begin{equation}
I(t,x)=-\frac{c}{|x|^{n}}\int_0^\ity \partial_r \Big(X^m \Big( e^{-c_1 r^{2\delta}t} \sin (c_2 r^\sigma t)\big(1-\chi(r)\big) r^{2\beta-\sigma+2m}\Big)\Big) \sin(r|x|)dr. \label{l3.1.2}
\end{equation}
A standard calculation leads to
\begin{align*}
I(t,x)&= \sum_{j=0}^m \sum_{k=0}^{j+1}\frac{c_{jk}}{|x|^{n}}\int_0^\ity \partial_r^{j+1-k} e^{-c_1 r^{2\delta}t}\, \partial_r^{k}\Big(\sin (c_2 r^\sigma t)\big(1-\chi(r)\big)\Big) r^{2\beta-\sigma+j} \sin(r|x|)dr\\
&\quad + \sum_{j=0}^m \sum_{k=0}^j \frac{c_{jk}}{|x|^{n}}\int_0^\ity \partial_r^{j-k} e^{-c_1 r^{2\delta}t}\, \partial_r^{k+1}\Big(\sin (c_2 r^\sigma t)\big(1-\chi(r)\big)\Big) r^{2\beta-\sigma+j} \sin(r|x|)dr\\
&\quad + \sum_{j=1}^m \sum_{k=0}^j \frac{c_{jk}}{|x|^{n}}\int_0^\ity \partial_r^{j-k} e^{-c_1 r^{2\delta}t}\, \partial_r^{k}\Big(\sin (c_2 r^\sigma t)\big(1-\chi(r)\big)\Big) r^{2\beta-\sigma+j-1} \sin(r|x|)dr
\end{align*}
with some constants $c_{jk}$. For this reason, we only need to study the integrals
\begin{equation}
I_{j,k}(t,x):= \int_0^\ity \partial_r^{j+1-k} e^{-c_1 r^{2\delta}t}\, \partial_r^{k}\Big(\sin (c_2 r^\sigma t)\big(1-\chi(r)\big)\Big) r^{2\beta-\sigma+j} \sin(r|x|)dr. \label{l3.1.3}
\end{equation}
Due to the large values of $r$, we can see that on the support of $1-\chi(r)$ and on the support of its derivatives it holds
\begin{align*}
\big|\partial_r^l e^{-c_1 r^{2\delta}t}\big| &\lesssim e^{-c_1 r^{2\delta}t}r^{l(2\delta-1)} t^{l},\\
\Big|\partial_r^{l}\Big(\sin (c_2 r^\sigma t)\big(1-\chi(r)\big)\Big)\Big| &\lesssim r^{l(\sigma-1)} t^{l} \text{ for } l= 0,\cdots,m.
\end{align*}
Hence, we imply for large $r$, $j= 0,\cdots,m$ and $k= 0,\cdots,j$ the estimates
$$\Big|\partial_r^{j+1-k} e^{-c_1 r^{2\delta}t}\, \partial_r^{k}\Big(\sin (c_2 r^\sigma t)\big(1-\chi(r)\big)\Big) r^{2\beta-\sigma+j}\Big| \lesssim e^{-c_1 r^{2\delta}t} t^{j+1} r^{2\delta (j+1)+k(\sigma-2\delta)+2\beta-1} $$
on the support of $1-\chi(r)$ and on the support of its derivatives. By splitting of the integral (\ref{l3.1.3})  into two parts, we get on the one hand
\begin{equation}
\Big |\int_0^{\frac{\pi}{2|x|}} \partial_r^{j+1-k} e^{-c_1 r^{2\delta}t}\, \partial_r^{k}\Big(\sin (c_2 r^\sigma t)\big(1-\chi(r)\big)\Big) r^{2\beta-\sigma+j} \sin(r|x|)dr \Big| \lesssim \frac{1}{|x|^{2\delta}}e^{-ct} \label{l3.1.4}
\end{equation}
for some constant $c>0$. On the other hand, we can carry out one more step of partial integration in the remaining integral as follows:
\begin{align}
&\Big|\int_{\frac{\pi}{2|x|}}^\ity \partial_r^{j+1-k} e^{-c_1 r^{2\delta}t}\, \partial_r^{k}\Big(\sin (c_2 r^\sigma t)\big(1-\chi(r)\big)\Big) r^{2\beta-\sigma+j} \sin(r|x|)dr \Big| \nonumber \\
&\quad \lesssim \frac{1}{|x|}\Big|\partial_r^{j+1-k} e^{-c_1 r^{2\delta}t}\, \partial_r^{k}\Big(\sin (c_2 r^\sigma t)\big(1-\chi(r)\big)\Big) r^{2\beta-\sigma+j} \cos(r|x|) \Big|_{r=\frac{\pi}{2|x|}}^\ity \nonumber \\
&\qquad + \frac{1}{|x|}\int_{\frac{\pi}{2|x|}}^\ity \Big|\partial_r \Big( \partial_r^{j+1-k} e^{-c_1 r^{2\delta}t}\, \partial_r^{k}\Big(\sin (c_2 r^\sigma t)\big(1-\chi(r)\big)\Big) r^{2\beta-\sigma+j}\Big) \cos(r|x|)\Big| \,dr \lesssim \frac{1}{|x|}e^{-ct} \label{l3.1.5}
\end{align}
for some constant $c>0$. Here we also note that for all $j= 0,\cdots,m$ and $k= 0,\cdots,j$ we have the estimates
$$\Big|\partial_r \Big( \partial_r^{j+1-k} e^{-c_1 r^{2\delta}t} \partial_r^{k}\Big(\sin (c_2 r^\sigma t)\big(1-\chi(r)\big)\Big) r^{2\beta-\sigma+j}\Big)\Big| \lesssim e^{-c_1 r^{2\delta}t} t^{j+2} r^{2\delta (j+1)+k(\sigma-2\delta)+2\beta-2}. $$
Therefore, from (\ref{l3.1.2}) to (\ref{l3.1.5}) we have produced terms $|x|^{-(n+2\delta)}$ and $|x|^{-(n+1)}$ which guarantee the $L^1$ property in $x$. Summarizing, implies for all $t \in [1,\ity)$ and $n=2m+1$ the estimates:
$$\Big\|F^{-1}\Big(e^{-c_1 |\xi|^{2\delta}t}|\xi|^{2\beta-\sigma} \sin (c_2 |\xi|^\sigma t) \big(1-\chi(|\xi|)\big)\Big)(t,\cdot)\Big\|_{L^1(|x|\ge 1)} \lesssim e^{-ct} \text{ for some }c>0. $$
\noindent Let us consider even spatial dimension $n=2m, m \ge 1$ in the first case $t\in [1,\ity)$. Carrying out $m-1$ steps of partial integration we obtain
\begin{align}
&I(t,x)= \frac{c}{|x|^{2m-2}}\int_0^\ity X^{m-1}\Big( e^{-c_1 r^{2\delta}t} \sin (c_2 r^\sigma t)\big(1-\chi(r)\big) r^{2\beta-\sigma+2m-1} \Big) \tilde{J}_0(r|x|) dr. \nonumber \\
&\qquad =\sum_{j=0}^{m-1}\frac{c_j}{|x|^{2m-2}}\int_0^\ity \partial_r^j \Big( e^{-c_1 r^{2\delta}t} \sin (c_2 r^\sigma t)\big(1-\chi(r)\big) r^{2\beta-\sigma} \Big) r^{j+1} \tilde{J}_0(r|x|) dr =:\sum_{j=0}^{m-1} c_j I_j(t,x). \label{l3.1.6}
\end{align}
Applying the first rule of modified Bessel functions for $\mu=1$ and the fifth rule for $\mu=0$ from Proposition \ref{PropertiesModifiedBesselfunctions}, after two more steps of partial integration we have
\begin{equation}
I_0(t,x)= -\frac{1}{|x|^{n}}\int_1^\ity \partial_r \Big( \partial_r \Big( e^{-c_1 r^{2\delta}t} \sin (c_2 r^\sigma t)\big(1-\chi(r)\big) r^{2\beta-\sigma}\Big) r \Big) \tilde{J}_0(r|x|) dr. \label{l3.1.7}
\end{equation}
Noting that for large $r$ and all $j=0,\cdots,m$ we have the inequality
$$ \Big|\partial^j_r \Big( e^{-c_1 r^{2\delta}t} \sin (c_2 r^\sigma t)\big(1-\chi(r)\big) r^{2\beta-\sigma}\Big)\Big| \lesssim e^{-c_1 r^{2\delta}t}t^j r^{j(\sigma-1)+2\beta-\sigma}. $$
Hence, we get
$$\Big|\partial_r \Big( \partial_r \Big( e^{-c_1 r^{2\delta}t} \sin (c_2 r^\sigma t)\big(1-\chi(r)\big) r^{2\beta-\sigma}\Big) r \Big)\Big| \lesssim e^{-c_1 r^{2\delta}t}t^2 r^{\sigma+2\beta-1}$$
on the support of $1-\chi(r)$. Now using the estimate $|\tilde{J}_0(s)| \le Cs^{-\frac{1}{2}}$ for $s>1$ we conclude
\begin{align}
&\Big| \int_1^\ity \partial_r \Big( \partial_r \Big( e^{-c_1 r^{2\delta}t} \sin (c_2 r^\sigma t)\big(1-\chi(r)\big) r^{2\beta-\sigma}\Big) r \Big) \tilde{J}_0(r|x|) dr \Big| \nonumber \\
&\quad \lesssim \int_1^\ity e^{-c_1 r^{2\delta}t}t^2 r^{\sigma+2\beta-1} \frac{1}{(r|x|)^{\frac{1}{2}}}dr = \frac{1}{|x|^{\frac{1}{2}}} t^2 \int_1^\ity e^{-c_1 r^{2\delta}t} r^{\sigma+2\beta-\frac{3}{2}} dr \lesssim \frac{1}{|x|^{\frac{1}{2}}}e^{-ct} \label{l3.1.8}
\end{align}
for some constant $c>0$. Therefore, from (\ref{l3.1.7}) and (\ref{l3.1.8}) we have $$\|I_0(t,\cdot )\|_{L^1(|x|\ge 1)}\lesssim e^{-ct} \text{ for all } t\in [1,\ity), \text{ and some constant }c>0. $$
Let $j \in [1,m-1]$ be an integer. By using again the first rule of modified Bessel functions for $\mu=1$ and the fifth rule for $\mu=0$ from Proposition \ref{PropertiesModifiedBesselfunctions} and carrying out partial integration we can re-write $I_j(t,x)$ in (\ref{l3.1.6}) as follows:
\begin{align*}
I_j(t,x)&= -\frac{1}{|x|^{2m}}\int_0^\ity \partial_r \Big( \partial_r^{j+1}\Big( e^{-c_1 r^{2\delta}t} \sin (c_2 r^\sigma t)\big(1-\chi(r)\big) r^{2\beta-\sigma}\Big) r^{j+1}\Big) \tilde{J}_0(r|x|) dr\\
&\quad -\frac{j}{|x|^{2m}}\int_0^\ity \partial_r \Big( \partial_r^j \Big( e^{-c_1 r^{2\delta}t} \sin (c_2 r^\sigma t)\big(1-\chi(r)\big) r^{2\beta-\sigma}\Big) r^j \Big) \tilde{J}_0(r|x|) dr.
\end{align*}
Applying an analogous treatment as we did for $I_0=I_0(t,x)$ implies
$$\|I_j(t,\cdot )\|_{L^1(|x|\ge 1)}\lesssim e^{-ct} \text{ for all } t\in [1,\ity) \text{ and } j=1,\cdots,m-1, $$
where $c$ is a suitable positive constant. Therefore, we have the following desired estimate for all $t\in [1,\ity)$ and $n=2m$:
$$\Big\| F^{-1}\Big(e^{-c_1 |\xi|^{2\delta}t}|\xi|^{2\beta-\sigma} \sin (c_2 |\xi|^\sigma t) \big(1-\chi(|\xi|)\big)\Big)(t,\cdot)\Big\|_{L^1(|x|\ge 1)}\lesssim  e^{-ct} \text{ for some constant }c>0. $$
\noindent Let us turn to the second case $t \in (0,1]$. By the change of variables $\xi=t^{-\frac{1}{2\delta}}\eta $ we get
\begin{align*}
&F^{-1}\Big(e^{-c_1 |\xi|^{2\delta}t}|\xi|^{2\beta-\sigma} \sin (c_2 |\xi|^\sigma t) \big(1-\chi(|\xi|)\big)\Big)(t,x)\\
&\qquad = t^{-\frac{n+2\beta-\sigma}{2\delta}}F^{-1}\Big(e^{-c_1 |\eta |^{2\delta}}|\eta |^{2\beta-\sigma} \sin (c_2 |\eta |^\sigma t^{1-\frac{\sigma}{2\delta}}) \big(1-\chi(t^{-\frac{1}{2\delta}}|\eta |)\big)\Big)(t,t^{-\frac{1}{2\delta}}x).
\end{align*}
Hence, we have
\begin{align*}
&\Big\|F^{-1}\Big(e^{-c_1 |\xi|^{2\delta}t}|\xi|^{2\beta-\sigma} \sin (c_2 |\xi|^\sigma t) \big(1-\chi(|\xi|)\big)\Big)(t,\cdot)\Big\|_{L^1}\\
&\qquad = t^{\frac{\sigma- 2\beta}{2\delta}}\Big\|F^{-1}\Big(e^{-c_1 |\eta |^{2\delta}}|\eta |^{2\beta-\sigma} \sin (c_2 |\eta |^\sigma t^{1-\frac{\sigma}{2\delta}}) \big(1-\chi(t^{-\frac{1}{2\delta}}|\eta |)\big)\Big)(t,\cdot)\Big\|_{L^1}.
\end{align*}
For this reason, we only need to study the Fourier multiplier in the form
$$H(t,x):= F^{-1}\Big(e^{-c_1 |\eta |^{2\delta}}|\eta |^{2\beta-\sigma} \sin (c_2 |\eta |^\sigma t^{1-\frac{\sigma}{2\delta}}) \big(1-\chi(t^{-\frac{1}{2\delta}}|\eta |)\big)\Big)(t,x). $$
First, we localize to small $|x|\le 1$. Then, we derive immediately
$$\|H(t,\cdot)\|_{L^1(|x| \le 1)} \lesssim t^{1-\frac{\sigma}{2\delta}}. $$
Therefore, we may conclude for $|x|\le 1$ the estimate
$$ \Big\|F^{-1}\Big(e^{-c_1 |\xi|^{2\delta}t}|\xi|^{2\beta-\sigma} \sin (c_2 |\xi|^\sigma t) \big(1-\chi(|\xi|)\big)\Big)(t,\cdot)\Big\|_{L^1(|x| \le 1)} \lesssim t^{1-\frac{\beta}{\delta}}. $$
We assume now $|x| \ge 1$. Using modified Bessel functions we shall estimate
\begin{equation}
H(t,x)=c \int_0^\ity e^{-c_1 r^{2\delta}}r^{2\beta-\sigma} \sin (c_2 r^\sigma t^{1-\frac{\sigma}{2\delta}}) \big(1-\chi(t^{-\frac{1}{2\delta}}r)\big) r^{n-1}\tilde{J}_{\frac{n}{2}-1}(r|x|) dr. \label{l3.1.9}
\end{equation}
Let us consider odd spatial dimensions $n=2m+1, m \ge 1.$ Then, carrying out $m+1$ steps of partial integration we re-write (\ref{l3.1.9}) as follows:
\begin{align}
&H(t,x)= \frac{c}{|x|^n} \int_0^\ity \partial_r \Big(X^m \Big( e^{-c_1 r^{2\delta}}\sin (c_2 r^\sigma t^{1-\frac{\sigma}{2\delta}}) \big(1-\chi(t^{-\frac{1}{2\delta}}r)\big) r^{2\beta-\sigma+2m}\Big) \Big) \sin(r|x|)dr \nonumber \\
&=:\sum_{1\le j+k \le m+1,\,\, j,\,k \ge 0} \frac{c_{jk}}{|x|^n} \int_0^\ity \partial_r^{j+1} e^{-c_1 r^{2\delta}} \partial^k_r \Big( \sin (c_2 r^\sigma t^{1-\frac{\sigma}{2\delta}}) \big(1-\chi(t^{-\frac{1}{2\delta}}r)\big) \Big) r^{j+k+2\beta-\sigma} \sin(r|x|)dr. \label{l3.1.10}
\end{align}
In order to estimate the function $H(t,x)$ we use the following auxiliary estimates:
\begin{align*}
\big|\partial_r^j e^{-c_1 r^{2\delta}}\big| &\lesssim
\begin{cases}
e^{-c_1 r^{2\delta}} \text{ if } j =0, &\\
e^{-c_1 r^{2\delta}}(r^{2\delta-j}+r^{j(2\delta-1)}) \lesssim e^{-c_1 r^{2\delta}}r^{2\delta-j}(1+r^{2\delta})^{j-1} \text{ if } j=1,\cdots,m,&
\end{cases}\\
\big|\partial_r^j \sin (c_2 r^\sigma t^{1-\frac{\sigma}{2\delta}})\big| &\lesssim
\begin{cases}
r^\sigma t^{1-\frac{\sigma}{2\delta}} \text{ if } j =0, &\\
r^{\sigma-j} t^{1-\frac{\sigma}{2\delta}}+(r^{\sigma-1} t^{1-\frac{\sigma}{2\delta}})^j \lesssim r^{\sigma-j} t^{1-\frac{\sigma}{2\delta}}(1+r^\sigma t^{1-\frac{\sigma}{2\delta}})^{j-1} \text{ if } j=1,\cdots,m.&
\end{cases}
\end{align*}
From the above estimates we may derive
$$\Big|\partial^k_r \Big( \sin (c_2 r^\sigma t^{1-\frac{\sigma}{2\delta}}) \big(1-\chi(t^{-\frac{1}{2\delta}}r)\big) \Big)\Big| \lesssim  \begin{cases}
r^\sigma t^{1-\frac{\sigma}{2\delta}} \text{ if } k =0, &\\
r^{\sigma-k} t^{1-\frac{\sigma}{2\delta}}(1+r^\sigma t^{1-\frac{\sigma}{2\delta}})^{k-1} \text{ if } k=1,\cdots,m.&
\end{cases} $$
Hence, we have
\begin{align*}
&\Big|\partial_r^{j+1} e^{-c_1 r^{2\delta}}\partial^k_r \Big( \sin (c_2 r^\sigma t^{1-\frac{\sigma}{2\delta}}) \big(1-\chi(t^{-\frac{1}{2\delta}}r)\big) \Big) r^{j+k+2\beta-\sigma}\Big| \\
&\qquad \lesssim
\begin{cases}
e^{-c_1 r^{2\delta}} t^{1-\frac{\sigma}{2\delta}} r^{2\delta+2\beta-1}(1+r^{2\delta})^j \text{ if } k =0, &\\
e^{-c_1 r^{2\delta}} t^{1-\frac{\sigma}{2\delta}} r^{2\delta+2\beta-1}(1+r^\sigma t^{1-\frac{\sigma}{2\delta}})^{k+j-1} \text{ if } k=1,\cdots,m,&
\end{cases}
\end{align*}
where we also note that $|\xi| \in [1,\ity)$, that is, $r \in [t^{\frac{1}{2\delta}},\ity)$ and $rt^{-\frac{1}{2\delta}}\ge 1$. Now, let us devote to $k=0$. By splitting the integral in (\ref{l3.1.10}) into two parts, on the one hand we obtain the following estimate for $t^{\frac{1}{2\delta}} < \frac{\pi}{2|x|}$:
\begin{align}
&\Big |\int_{t^{\frac{1}{2\delta}}}^{\frac{\pi}{2|x|}} \partial_r^{j+1} e^{-c_1 r^{2\delta}}\sin (c_2 r^\sigma t^{1-\frac{\sigma}{2\delta}}) \big(1-\chi(t^{-\frac{1}{2\delta}}r)\big) r^{j+2\beta-\sigma} \sin(r|x|)dr \Big| \nonumber \\
&\qquad \lesssim  t^{1-\frac{\sigma}{2\delta}} \int_{t^{\frac{1}{2\delta}}}^{\frac{\pi}{2|x|}} r^{2\delta+2\beta-1}(1+r^{2\delta})^j dr \lesssim t^{1-\frac{\sigma}{2\delta}} \Big(\frac{1}{|x|^{2\delta+2\beta}}+ \frac{1}{|x|^{2\delta(j+1)+2\beta}} \Big) \lesssim \frac{t^{1-\frac{\sigma}{2\delta}}}{|x|^{2\delta+2\beta}}. \label{l3.1.11}
\end{align}
On the other hand, carrying out one more step of partial integration we derive
\begin{align}
&\Big|\int_{\frac{\pi}{2|x|}}^\ity \partial_r^{j+1} e^{-c_1 r^{2\delta}}\sin (c_2 r^\sigma t^{1-\frac{\sigma}{2\delta}}) \big(1-\chi(t^{-\frac{1}{2\delta}}r)\big) r^{j+2\beta-\sigma} \sin(r|x|) dr \Big| \nonumber \\
&\quad \lesssim \frac{1}{|x|}\left |\partial_r^{j+1} e^{-c_1 r^{2\delta}}\sin (c_2 r^\sigma t^{1-\frac{\sigma}{2\delta}}) \big(1-\chi(t^{-\frac{1}{2\delta}}r)\big) r^{j+2\beta-\sigma} \cos(r|x|) \right |_{r=\frac{\pi}{2|x|}}^\ity \nonumber \\
&\qquad + \frac{1}{|x|}\int_{\frac{\pi}{2|x|}}^\ity \Big|\partial_r \Big(\partial_r^{j+1} e^{-c_1 r^{2\delta}}\sin (c_2 r^\sigma t^{1-\frac{\sigma}{2\delta}}) \big(1-\chi(t^{-\frac{1}{2\delta}}r)\big) r^{j+2\beta-\sigma}\Big) \cos(r|x|)\Big| dr \nonumber \\
&\quad \lesssim \frac{t^{1-\frac{\sigma}{2\delta}}}{|x|}\int_{\frac{\pi}{2|x|}}^\ity e^{-c_1 r^{2\delta}} r^{2\delta+2\beta-2}(1+r^{2\delta})^{j+1} dr \lesssim \frac{t^{1-\frac{\sigma}{2\delta}}}{|x|} \Big( \int_{\frac{\pi}{2|x|}}^1 r^{2\delta+2\beta-2} dr + \int_1^\ity e^{-c_1 r^{2\delta}} r^{2\delta(j+2)+2\beta-2} dr \Big) \nonumber \\
&\quad \lesssim
\begin{cases}
t^{1-\frac{\sigma}{2\delta}}\Big(\frac{1}{|x|}+ \frac{1}{|x|^{2\delta+2\beta}}\Big) \text{ if }2\delta+2\beta \neq 1&\\
\frac{t^{1-\frac{\sigma}{2\delta}}}{|x|}\log(e+|x|) \text{ if }2\delta+2\beta= 1 &
\end{cases} \lesssim t^{1-\frac{\sigma}{2\delta}}\Big(\frac{1}{|x|^{\frac{1}{2}}}+ \frac{1}{|x|^{2\delta+2\beta}}\Big), \label{l3.1.12}
\end{align}
where we also note that
$$\Big|\partial_r \Big( \partial_r^{j+1} e^{-c_1 r^{2\delta}}\sin (c_2 r^\sigma t^{1-\frac{\sigma}{2\delta}}) \big(1-\chi(t^{-\frac{1}{2\delta}}r)\big) r^{j+2\beta-\sigma}\Big)\Big| \lesssim e^{-c_1 r^{2\delta}} t^{1-\frac{\sigma}{2\delta}} r^{2\delta+2\beta-2}(1+r^{2\delta})^{j+1}. $$
For $k=1,\cdots,m$, after an analogous treatment as we did for $k=0$ we get
\begin{equation}
\Big |\int_{t^{\frac{1}{2\delta}}}^{\frac{\pi}{2|x|}} \partial_r^{j+1} e^{-c_1 r^{2\delta}}\partial^k_r \Big( \sin (c_2 r^\sigma t^{1-\frac{\sigma}{2\delta}}) \big(1-\chi(t^{-\frac{1}{2\delta}}r)\big) \Big) r^{j+k+2\beta-\sigma} \sin(r|x|)dr \Big| \lesssim \frac{t^{(k+j)(1-\frac{\sigma}{2\delta})}}{|x|^{2\delta+2\beta}}, \label{l3.1.13}
\end{equation}
and
\begin{eqnarray}
&& \Big|\int_{\frac{\pi}{2|x|}}^\ity \partial_r^{j+1} e^{-c_1 r^{2\delta}}\partial^k_r \Big( \sin (c_2 r^\sigma t^{1-\frac{\sigma}{2\delta}}) \big(1-\chi(t^{-\frac{1}{2\delta}}r)\big) \Big) r^{j+k+2\beta-\sigma} \sin(r|x|) dr \Big| \nonumber \\ && \qquad \lesssim t^{(k+j+1)(1-\frac{\sigma}{2\delta})}\Big(\frac{1}{|x|^{\frac{1}{2}}}+ \frac{1}{|x|^{2\delta+2\beta}}\Big), \label{l3.1.14}
\end{eqnarray}
where we can see that
$$ \Big|\partial_r \Big( \partial_r^{j+1} e^{-c_1 r^{2\delta}}\partial^k_r \Big( \sin (c_2 r^\sigma t^{1-\frac{\sigma}{2\delta}}) \big(1-\chi(t^{-\frac{1}{2\delta}}r)\big) \Big) r^{j+k+2\beta-\sigma}\Big)\Big| \lesssim e^{-c_1 r^{2\delta}} t^{1-\frac{\sigma}{2\delta}} r^{2\delta+2\beta-2}(1+r^\sigma t^{1-\frac{\sigma}{2\delta}})^{k+j}. $$
Hence, from (\ref{l3.1.10}) to (\ref{l3.1.14}) we have produced terms $|x|^{-(n+2\delta+2\beta)}$ and $|x|^{-(n+\frac{1}{2})}$ which guarantee the $L^1$ property in $x$. For this reason we arrive at for all $t \in (0,1]$ and $n=2m+1$ the following estimates:
$$\Big\| F^{-1}\Big(e^{-c_1 |\xi|^{2\delta}t}|\xi|^{2\beta-\sigma} \sin (c_2 |\xi|^\sigma t) \big(1- \chi(|\xi|)\big) \Big)(t,\cdot)\Big\|_{L^1(|x|\ge 1)} \lesssim t^{-(m+2)(\frac{\sigma}{2\delta}-1)+\frac{\sigma-2\beta}{2\delta}}. $$
Let us consider even spatial dimensions $n=2m, m \ge 1$. Carrying out $m-1$ steps of partial integration we re-write (\ref{l3.1.9}) as follows:
\begin{align*}
&H(t,x)= \frac{c}{|x|^{2m-2}} \int_0^\ity X^{m-1} \Big(e^{-c_1 r^{2\delta}}\sin (c_2 r^\sigma t^{1-\frac{\sigma}{2\delta}}) \big(1-\chi(t^{-\frac{1}{2\delta}}r)\big) r^{2\beta-\sigma+2m-1} \Big) \tilde{J}_0 (r|x|)dr \\
&\quad =\sum_{j=0}^{m-1}\frac{c_j}{|x|^{2m-2}} \int_0^\ity \partial_r^j \Big( e^{-c_1 r^{2\delta}}\sin (c_2 r^\sigma t^{1-\frac{\sigma}{2\delta}}) \big(1-\chi(t^{-\frac{1}{2\delta}}r)\big) r^{2\beta-\sigma}\Big)  r^{j+1} \tilde{J}_0 (r|x|)dr =:\sum_{j=0}^{m-1} c_j I_j(t,x).
\end{align*}
Using the first rule of modified Bessel functions for $\mu=1$ and the fifth rule for $\mu=0$ from Proposition \ref{PropertiesModifiedBesselfunctions} and performing two more steps of partial integration we get
$$|I_0(t,x)|= \frac{1}{|x|^{2m}} \int_0^\ity \Big|\partial_r \Big( \partial_r \Big( e^{-c_1 r^{2\delta}}\sin (c_2 r^\sigma t^{1-\frac{\sigma}{2\delta}}) \big(1-\chi(t^{-\frac{1}{2\delta}}r)\big) r^{2\beta-\sigma}\Big) r \Big) \tilde{J}_0 (r|x|)\Big| dr. $$
We can see that for $j=1,\cdots,m$ we have
$$\Big| \partial^j_r \Big( e^{-c_1 r^{2\delta}}\sin (c_2 r^\sigma t^{1-\frac{\sigma}{2\delta}}) \big(1-\chi(t^{-\frac{1}{2\delta}}r)\big) r^{2\beta-\sigma}\Big)\Big| \lesssim e^{-c_1 r^{2\delta}} t^{1-\frac{\sigma}{2\delta}}(1+ r^\sigma t^{1-\frac{\sigma}{2\delta}})^{j-1}(r^{2\delta+2\beta-j}+ r^{2\beta-j}) $$
on the support of $1-\chi(t^{-\frac{1}{2\delta}}r)$ and on the support of its derivatives. Therefore, we may conclude
$$\Big|\partial_r \Big( \partial_r \Big( e^{-c_1 r^{2\delta}}\sin (c_2 r^\sigma t^{1-\frac{\sigma}{2\delta}}) \big(1-\chi(t^{-\frac{1}{2\delta}}r)\big) r^{2\beta-\sigma}\Big) r \Big)\Big| \lesssim e^{-c_1 r^{2\delta}} t^{1-\frac{\sigma}{2\delta}}(1+ r^\sigma t^{1-\frac{\sigma}{2\delta}}) (r^{2\delta+2\beta-1}+r^{2\beta-1}). $$
Since $|\tilde{J}_0(s)| \le C$ for $s \in [0,1]$ we obtain for $t^{\frac{1}{2\delta}}<\frac{1}{|x|}$ the estimate
\begin{align}
&\int_{t^{\frac{1}{2\delta}}}^{\frac{1}{|x|}} \Big|\partial_r \Big( \partial_r \Big( e^{-c_1 r^{2\delta}}\sin (c_2 r^\sigma t^{1-\frac{\sigma}{2\delta}}) \big(1-\chi(t^{-\frac{1}{2\delta}}r)\big) r^{2\beta-\sigma}\Big) r \Big) \tilde{J}_0(r|x|)\Big| dr \nonumber \\
&\qquad \lesssim  t^{1-\frac{\sigma}{2\delta}} \int_{t^{\frac{1}{2\delta}}}^{\frac{1}{|x|}} e^{-c_1 r^{2\delta}}(1+ r^\sigma t^{1-\frac{\sigma}{2\delta}}) r^{2\beta-1}dr \,\,\, \Big(\text{since }  r^{2\delta} \le 1 \text{ for } r \le \frac{1}{|x|} \Big) \nonumber \\
& \qquad\lesssim  t^{1-\frac{\sigma}{2\delta}} \int_{t^{\frac{1}{2\delta}}}^{\frac{1}{|x|}}  r^{2\beta-1}dr+ t^{2(1-\frac{\sigma}{2\delta})} \int_{t^{\frac{1}{2\delta}}}^{\frac{1}{|x|}}r^{\sigma+ 2\beta-1}dr \nonumber \\
&\qquad\lesssim  t^{2(1-\frac{\sigma}{2\delta})} \int_{t^{\frac{1}{2\delta}}}^{\frac{1}{|x|}}  r^{2\beta-1+\sigma-2\delta} dr+ \frac{t^{2(1-\frac{\sigma}{2\delta})}}{|x|^{\sigma+2\beta}} \,\,\, \Big(\text{since} \,\,r^{2\delta-\sigma} \le t^{1-\frac{\sigma}{2\delta}} \text{ for } r \ge t^{\frac{1}{2\delta}} \Big) \nonumber \\
& \qquad\lesssim  t^{2(1-\frac{\sigma}{2\delta})} \Big(\frac{1}{|x|^{\sigma- 2\delta+2\beta}}+ \frac{1}{|x|^{\sigma+2\beta}}\Big) \lesssim  t^{2(1-\frac{\sigma}{2\delta})} \frac{1}{|x|^{\sigma- 2\delta+2\beta}}. \label{l3.1.16}
\end{align}
Moreover, we use $|\tilde{J}_0(s)| \le Cs^{-\frac{1}{2}}$ for $s>1$ to conclude
\begin{align}
&\int_{\frac{1}{|x|}}^\ity \Big|\partial_r \Big( \partial_r \Big( e^{-c_1 r^{2\delta}}\sin (c_2 r^\sigma t^{1-\frac{\sigma}{2\delta}}) \big(1-\chi(t^{-\frac{1}{2\delta}}r)\big) r^{2\beta-\sigma}\Big) r \Big) \tilde{J}_0(r|x|)\Big| dr \nonumber \\
&\qquad\lesssim \frac{t^{1-\frac{\sigma}{2\delta}}}{|x|^{\frac{1}{2}}} \int_{\frac{1}{|x|}}^\ity e^{-c_1 r^{2\delta}} (1+ r^\sigma t^{1-\frac{\sigma}{2\delta}}) ( r^{2\delta+2\beta-\frac{3}{2}}+ r^{2\beta-\frac{3}{2}}) dr \nonumber \\
&\qquad\lesssim \frac{t^{1-\frac{\sigma}{2\delta}}}{|x|^{\frac{1}{2}}} \int_{\frac{1}{|x|}}^1 r^{2\beta-\frac{3}{2}}dr+ \frac{t^{2(1-\frac{\sigma}{2\delta})}}{|x|^{\frac{1}{2}}} \lesssim \frac{t^{2(1-\frac{\sigma}{2\delta})}}{|x|^{\frac{1}{2}}} \int_{\frac{1}{|x|}}^1 r^{2\beta-\frac{3}{2}+\sigma-2\delta} dr+ \frac{t^{2(1-\frac{\sigma}{2\delta})}}{|x|^{\frac{1}{2}}} \,\,\, \Big(\text{since} \,\,r^{2\delta-\sigma} \le t^{1-\frac{\sigma}{2\delta}}\Big) \nonumber \\
&\qquad\lesssim
\begin{cases}
t^{2(1-\frac{\sigma}{2\delta})}\Big(\frac{1}{|x|^{\frac{1}{2}}}+ \frac{1}{|x|^{\sigma-2\delta+2\beta}}\Big) \text{ if }\sigma- 2\delta+ 2\beta \neq \frac{1}{2}&\\
\frac{t^{2(1-\frac{\sigma}{2\delta})}}{|x|^{\frac{1}{2}}}\log(e+|x|) \text{ if }\sigma- 2\delta+ 2\beta= \frac{1}{2}&
\end{cases} \lesssim t^{2(1-\frac{\sigma}{2\delta})} \Big(\frac{1}{|x|^{\frac{1}{4}}}+ \frac{1}{|x|^{\sigma-2\delta+2\beta}}\Big). \label{l3.1.17}
\end{align}
Hence, from (\ref{l3.1.16}) and (\ref{l3.1.17}) we have produced terms $|x|^{-(n+\frac{1}{4})}$ and $|x|^{-(n+\sigma-2\delta+2\beta)}$ which guarantee the $L^1$ property in $x$. Summarizing, we arrive at the estimate
$$\|I_0(t,\cdot)\|_{L^1(|x|\ge 1)} \lesssim  t^{2(1-\frac{\sigma}{2\delta})} \,\,\text{ for all } \,\,t\in (0,1]. $$
Let $j \in [1,m-1]$ be an integer. Then, repeating the above arguments we also derive for $t \in (0,1]$
$$ \|I_j(t,\cdot)\|_{L^1(|x|\ge 1)} \lesssim  t^{(j+2)(1-\frac{\sigma}{2\delta})}. $$
Therefore, we have proved that for all $t \in (0,1]$ and $n=2m$ the following estimates hold:
$$\Big\| F^{-1}\Big(e^{-c_1 |\xi|^{2\delta}t}|\xi|^{2\beta-\sigma} \sin (c_2 |\xi|^\sigma t) \big(1-\chi(|\xi|)\big)\Big)(t,\cdot) \Big\|_{L^1(|x|\ge 1)} \lesssim t^{-(m+1)(\frac{\sigma}{2\delta}-1)+\frac{\sigma-2\beta}{2\delta}}. $$
Summarizing, the proof to Lemma \ref{bd3.1.3} is completed.
\end{proof}

\begin{nx}
\fontshape{n}
\selectfont
In the proof of Lemma \ref{bd3.1.3} we explained our considerations for $n \ge 2$. Nevertheless, repeating the steps of the proof for odd spatial dimension we conclude that the statements of this lemma also hold for $n=1$. Here in the latter case we notice that we only carry out partial integration with no necessity to introduce the vector field $Xf(r)$ as we did in (\ref{l3.1.2}) and (\ref{l3.1.10}).
\end{nx}

\noindent Following the proof of Lemma \ref{bd3.1.3} we may conclude the following $L^1$ estimates, too.
\bbd \label{bd3.1.4}
The following estimates hold in $\R^n$ for any $n \ge 1$:
$$\Big\|F^{-1}\Big(e^{-c_1 |\xi|^{2\delta}t}|\xi|^{2\beta}\cos (c_2 |\xi|^\sigma t) \big(1-\chi(|\xi|)\big)\Big)(t,\cdot)\Big\|_{L^1}
\lesssim \begin{cases}
t^{-(2+[\frac{n}{2}])(\frac{\sigma}{2\delta}-1)-\frac{\beta}{\delta}} \text{ for } t\in (0,1], &\\
e^{-ct} \text{ for } t\in[1,\ity), &
\end{cases}$$
where $\beta \ge 0$ and $c$ is a suitable positive constant. Moreover, $c_1$ is a positive and $c_2 \ne 0$ is a real constant.
\ebd

\noindent Finally, we consider oscillating integrals with a more complicated oscillating integrand. We are going to prove the following result.
\bbd \label{bd3.1.5}
The following estimates hold in $\R^n$ for any $n \ge 1$:
$$\Big\|F^{-1}\Big(e^{-c_1 |\xi|^{2\delta}t}|\xi|^{2\beta} \frac{\sin \big(c_2 |\xi|^\sigma f(|\xi|) t \big)}{|\xi|^\sigma f(|\xi|)}\big(1-\chi(|\xi|)\big)\Big)(t,\cdot)\Big\|_{L^1} \lesssim \begin{cases}
t^{-(2+[\frac{n}{2}])(\frac{\sigma}{2\delta}-1)+\frac{\sigma-2\beta}{2\delta}} \text{ for } t\in (0,1], &\\
e^{-ct} \text{ for } t\in[1,\ity), &
\end{cases} $$
where \[ f(|\xi|)=\sqrt{1-\frac{1}{4|\xi|^{2\sigma-4\delta}}}, \] $\beta \ge 0$ and $c$ is a suitable positive constant. Moreover, $c_1$ is a positive and $c_2 \ne 0$ is a real constant.
\ebd
\begin{proof} The proof of this lemma is similar to the proof to Lemma \ref{bd3.1.3}. For this reason, we only present the steps which are different. Then, we shall repeat some of the arguments as we did in the proof to Lemma \ref{bd3.1.3} to conclude the desired estimates.\medskip

\noindent First, let us consider $|x| \ge 1$ and $t \in [1,\ity)$. In order to obtain exponential decay estimates in both cases of odd spatial dimensions $n=2m+1$ and even spatial dimensions $n=2m$ with $m \ge 1$, we shall prove the following estimates on the support of $1-\chi(r)$ and on the support of its derivatives:
$$ \Big| \partial^{k}_r \Big( \frac{\sin \big(c_2 r^\sigma f(r) t\big)}{f(r)}\big(1-\chi(r)\big) \Big) \Big| \lesssim r^{k(\sigma-1)}t^k \text{ for } k= 1,\cdots,m, $$
where
$$f(r)=\sqrt{1-\frac{1}{4r^{2\sigma-4\delta}}}. $$
Indeed, we shall apply Fa\`{a} di Bruno's formula as a main tool. We divide the proof of the above estimates into several sub-steps as follows: \medskip

\noindent Step 1: $\quad$ Applying Proposition \ref{FadiBruno'sformula1} with $h(s)= \sqrt{s}$ and $g(r)=1-\frac{1}{4}r^{-2(\sigma-2\delta)}$, we get
\begin{align*}
&\big|\partial_r^k f(r)\big| \lesssim \Big| \sum_{\substack{1\cdot m_1+\cdots+k\cdot m_k=k,\,m_i \ge 0}} g(r)^{\frac{1}{2}-(m_1+\cdots+m_k)}\prod_{j=1}^k \Big(-\frac{1}{4}r^{-2(\sigma-2\delta)-j}\Big)^{m_j}\Big|\\
&\qquad \lesssim \sum_{\substack{1\cdot m_1+\cdots+k\cdot m_k=k,\, m_i \ge 0}} r^{-2(\sigma-2\delta)(m_1+\cdots+m_k)-k}\lesssim r^{-k} \,\,\Big(\text{since}\,\,\frac{3}{4} \le g(r) \le 1 \,\,\,\text{for } \,r \ge 1 \Big).
\end{align*}
An analogous treatment gives
\begin{equation}
\Big|\partial_r^k \Big(\frac{1}{f(r)} \Big)\Big| \lesssim r^{-k} \text{ for }k=1,\cdots,m. \label{l3.3.1}
\end{equation}
\noindent Step 2: $\quad$ Applying Proposition \ref{FadiBruno'sformula1} with $h(s)= \sin(c_2\,s)$ and $g(r)= r^\sigma f(r) t$, we obtain
\begin{align}
&\big| \partial_r^k \sin \big(c_2 r^\sigma f(r) t\big) \big| \lesssim \Big| \sum_{\substack{1\cdot m_1+\cdots+k\cdot m_k=k,\, m_i \ge 0}} \sin \big(c_2 r^\sigma f(r) t\big)^{(m_1+\cdots+m_k)} \prod_{j=1}^k \Big( \partial_r^j \big(r^\sigma f(r) t\big) \Big)^{m_j} \Big| \nonumber \\
&\qquad \lesssim \Big| \sum_{\substack{1\cdot m_1+\cdots+k\cdot m_k=k,\, m_i \ge 0}} \,\, \prod_{j=1}^k \Big( t\, \sum_{l=0}^j C_j^l r^{\sigma-j+l} f(r)^{(l)} \Big)^{m_j} \Big| \nonumber \\
&\qquad \lesssim \sum_{\substack{1\cdot m_1+\cdots+k\cdot m_k=k,\, m_i \ge 0}} \,\, \prod_{j=1}^k (t\, r^{\sigma-j})^{m_j} \lesssim \sum_{\substack{1\cdot m_1+\cdots+k\cdot m_k=k,\, m_i \ge 0}} r^{-k} (t\, r^\sigma)^{m_1+\cdots+m_k} \lesssim t^k r^{k(\sigma-1)}. \label{l3.3.2}
\end{align}
Hence, from (\ref{l3.3.1}) and (\ref{l3.3.2}) using the product rule for higher derivatives we may conclude
$$\Big| \partial_r^k \Big(\frac{\sin \big(c_2 r^\sigma f(r) t\big)}{f(r)}\Big) \Big| \lesssim  t^k r^{k(\sigma-1)} \,\,\text{ for } \,\,k=1,\cdots,m. $$

\noindent Next, let us turn to the case $|x| \ge 1$ and $t \in (0,1]$. In order to prove the desired estimates by using similar ideas as in the proof to Lemma \ref{bd3.1.3}, we need to assert the following auxiliary estimates on the support of $1-\chi(t^{-\frac{1}{2\delta}}r)$ and on the support of its partial derivatives:
$$ \Big| \partial_r^k \Big( \frac{\sin \big(c_2 r^\sigma f(r) t^{1-\frac{\sigma}{2\delta}}\big)}{f(r)} \Big) \Big| \lesssim t^{1-\frac{\sigma}{2\delta}}r^{\sigma-k}(1+r^\sigma t^{1-\frac{\sigma}{2\delta}})^{k-1} \text{ if } k=1,\cdots,m, $$
where
$$f(r)= \sqrt{1-\frac{1}{4}t^{\frac{\sigma-2\delta}{\delta}}r^{-2(\sigma-2\delta)}}. $$
Indeed, we shall divide our proof into several sub-steps as follows:\medskip

\noindent Step 1: $\quad$ Applying Proposition \ref{FadiBruno'sformula1} with $h(s)=\sqrt{s}$ and $g(r)=1-\frac{1}{4}t^{\frac{\sigma-2\delta}{\delta}} r^{-2(\sigma-2\delta)}$, we obtain
\begin{align*}
&\big|\partial_r^k f(r)\big| \lesssim \Big| \sum_{\substack{1\cdot m_1+\cdots+k\cdot m_k=k,\, m_i \ge 0}} g(r)^{\frac{1}{2}-(m_1+\cdots+m_k)}\prod_{j=1}^k \Big(-\frac{1}{4}t^{\frac{\sigma-2\delta}{\delta}}r^{-2(\sigma-2\delta)-j}\Big)^{m_j}\Big|\\
&\qquad \lesssim \sum_{\substack{1\cdot m_1+\cdots+k\cdot m_k=k,\, m_i \ge 0}} \big( t^{\frac{\sigma-2\delta}{\delta}} r^{-2(\sigma-2\delta)}\big)^{m_1+\cdots+m_k} r^{-k} \,\,\,  \Big(\text{since } \frac{3}{4} \le g(r) \le 1 \text{ for } r \ge t^{\frac{1}{2\delta}}\Big)\\
&\qquad \lesssim r^{-k} \sum_{\substack{1\cdot m_1+\cdots+k\cdot m_k=k,\, m_i \ge 0}} (t^{-\frac{1}{2\delta}}r)^{-2(\sigma-2\delta)(m_1+\cdots+m_k)} \lesssim r^{-k} \,\,\,  \Big(\text{since } t^{-\frac{1}{2\delta}}r \ge 1 \text{ for } r \ge t^{\frac{1}{2\delta}}\Big).
\end{align*}
In an analogous way we may derive the estimates
\begin{equation}
\Big| \partial_r^k \Big(\frac{1}{f(r)} \Big)\Big|  \lesssim r^{-k} \text{ for } k=1,\cdots,m. \label{l3.3.3}
\end{equation}

\noindent Step 2: $\quad$ Applying Proposition \ref{FadiBruno'sformula1} with $h(s)=\sin(c_2\,s)$ and $g(r)=r^\sigma f(r) t^{1-\frac{\sigma}{2\delta}}$, we derive
\begin{align}
&\big| \partial_r^k \sin \big(c_2 r^\sigma f(r) t^{1-\frac{\sigma}{2\delta}}\big) \big| \nonumber \\
&\quad \lesssim \Big| \sum_{\substack{1\cdot m_1+\cdots+k\cdot m_k=k,\, m_i \ge 0}} \sin \big(c_2 r^\sigma f(r) t^{1-\frac{\sigma}{2\delta}}\big)^{(m_1+\cdots+m_k)} \,\, \prod_{j=1}^k \Big( \partial_r^j \big(r^\sigma f(r) t^{1-\frac{\sigma}{2\delta}}\big) \Big)^{m_j} \Big| \nonumber \\
&\quad \lesssim \Big| \sum_{\substack{1\cdot m_1+\cdots+k\cdot m_k=k,\, m_i \ge 0}} \,\, \prod_{j=1}^k \Big(t^{1-\frac{\sigma}{2\delta}} \sum_{l=0}^j C_j^l r^{\sigma-j+l} \partial_r^l f(r) \Big)^{m_j} \Big| \nonumber \\
&\quad \lesssim \sum_{\substack{1\cdot m_1+\cdots+k\cdot m_k=k,\, m_i \ge 0}} \,\, \prod_{j=1}^k (t^{1-\frac{\sigma}{2\delta}} r^{\sigma-j})^{m_j } \,\, \lesssim t^{1-\frac{\sigma}{2\delta}} r^{\sigma-k} \sum_{\substack{1\cdot m_1+\cdots+k\cdot m_k=k,\, m_i \ge 0}} (t^{1-\frac{\sigma}{2\delta}} r^\sigma)^{m_1+\cdots+m_k-1} \nonumber \\
&\quad \lesssim t^{1-\frac{\sigma}{2\delta}} r^{\sigma-k} \big( 1+(t^{1-\frac{\sigma}{2\delta}} r^\sigma)^{k-1}\big) \lesssim t^{1-\frac{\sigma}{2\delta}} r^{\sigma-k} \big( 1+t^{1-\frac{\sigma}{2\delta}} r^\sigma\big)^{k-1}. \label{l3.3.4}
\end{align}
Hence, from (\ref{l3.3.3}) and (\ref{l3.3.4}) using the product rule for higher derivatives we may conclude
$$\Big| \partial_r^k \Big(\frac{\sin \big(c_2 r^\sigma f(r) t^{1-\frac{\sigma}{2\delta}}\big)}{f(r)} \Big)\Big| \lesssim t^{1-\frac{\sigma}{2\delta}} r^{\sigma-k} \big( 1+t^{1-\frac{\sigma}{2\delta}} r^\sigma\big)^{k-1} \text{ for } k=1,\cdots,m. $$
Summarizing, the proof of Lemma \ref{bd3.1.5} is completed.
\end{proof}

\noindent Following the steps of the proof to Lemma \ref{bd3.1.5} we may prove the following statement, too.
\bbd \label{bd3.1.6}
The following estimates hold in $\R^n$ for any $n \ge 1$:
$$\Big\|F^{-1}\Big(e^{-c_1 |\xi|^{2\delta}t}|\xi|^{2\beta}\cos \big(c_2 |\xi|^\sigma f(|\xi|) t\big) \big(1-\chi(|\xi|)\big)\Big)(t,\cdot)\Big\|_{L^1}
\lesssim \begin{cases}
t^{-(2+[\frac{n}{2}])(\frac{\sigma}{2\delta}-1)-\frac{\beta}{\delta}} \text{ for } t\in (0,1],  &\\
e^{-ct} \text{ for } t\in[1,\ity), &
\end{cases}$$
where \[ f(|\xi|)=\sqrt{1-\frac{1}{4|\xi|^{2\sigma-4\delta}}}, \] $\beta \ge 0$ and $c$ is a suitable positive constant. Moreover, $c_1$ is a positive and $c_2 \ne 0$ is a real constant.
\ebd

\renewcommand{\proofname}{Proof of Proposition \ref{md3.1.4}.}
\begin{proof}
Thanks to the relation (\ref{p3.2.1}), to prove the first statement we choose $2\beta= a+2\delta$ and $2\beta= a$, respectively, in Lemmas \ref{bd3.1.5} and \ref{bd3.1.6}. Moreover, using the relation (\ref{p3.2.2}) and replacing $2\beta= a$ in Lemma \ref{bd3.1.5} we can conclude the second statement.
\end{proof}

\noindent Finally, from the statements of Propositions \ref{md3.1.2} and \ref{md3.1.4}, we may conclude the following $L^1$ estimates.
\bmd \label{md3.1.6}
The following estimates hold in $\R^n$ for any $n \ge 1$:
\begin{align*}
&\big\| F^{-1}\big(|\xi|^a \hat{K_0}\big)(t,\cdot)\big\|_{L^1}\lesssim \begin{cases}
t^{-(2+[\frac{n}{2}])(\frac{\sigma}{2\delta}-1)-\frac{a}{2\delta}} \text{ for } t\in (0,1], &\\
t^{-\frac{a}{2(\sigma-\delta)}} \text{ for } t\in[1,\ity), &
\end{cases}\\
&\big\| F^{-1}\big(|\xi|^a \hat{K_1}\big)(t,\cdot)\big\|_{L^1}\lesssim \begin{cases}
t^{1-(1+[\frac{n}{2}])(\frac{\sigma}{2\delta}-1)-\frac{a}{2\delta}} \text{ for } t\in (0,1], &\\
t^{1-\frac{a}{2(\sigma-\delta)}} \text{ for } t\in[1,\ity), &
\end{cases}
\end{align*}
for any non-negative number $a$.
\emd

\subsection{$L^\ity$ estimates}
\bmd \label{md3.1.9}
The following estimates hold in $\R^n$ for any $n \ge 1$:
\begin{align*}
&\big\|F^{-1}\big(|\xi|^a \hat{K_0}\chi(|\xi|)\big)(t,\cdot)\big\|_{L^\ity}\lesssim
\begin{cases}
1 \text{ for } t \in (0,1], &\\
 t^{-\frac{n+a}{2(\sigma-\delta)}} \text{ for } t \in [1,\ity), &
\end{cases}\\
&\Big\|F^{-1}\Big(|\xi|^a \hat{K_0}\big(1-\chi(|\xi|)\big)\Big)(t,\cdot)\Big\|_{L^\ity}\lesssim t^{-\frac{n+a}{2\delta}} \text{ for } t \in (0,\ity),\\
&\big\|F^{-1}\big(|\xi|^a \hat{K_1}\chi(|\xi|)\big)(t,\cdot)\big\|_{L^\ity}\lesssim
\begin{cases}
t \text{ for } t \in (0,1], &\\
t^{1-\frac{n+a}{2(\sigma-\delta)}} \text{ for } t \in [1,\ity), &
\end{cases}\\
&\Big\|F^{-1}\Big(|\xi|^a \hat{K_1}\big(1-\chi(|\xi|)\big)\Big)(t,\cdot)\Big\|_{L^\ity}\lesssim t^{1-\frac{n+a}{2\delta}} \text{ for } t \in (0,\ity),
\end{align*}
for any non-negative number $a$.
\emd

\renewcommand{\proofname}{Proof.}
\begin{proof}
First, taking account of the representation for $\hat{K_1}$ we can re-write it as follows:
$$ \hat{K_1}(t,\xi)= e^{\lambda_1 t}\f{1- e^{(\lambda_2- \lambda_1)t}}{\lambda_1-\lambda_2}=
\begin{cases}
te^{\lambda_1 t}\int_0^1 e^{-\theta \sqrt{|\xi|^{4\delta}-4|\xi|^{2\sigma}}t}d\theta \text{ for small } |\xi|, &\\
te^{\lambda_1 t}\int_0^1 e^{-\theta i \sqrt{4|\xi|^{2\sigma}-|\xi|^{4\delta}}t}d\theta \text{ for large } |\xi|. &
\end{cases} $$
For the sake of the asymptotic behavior of the characteristic roots in (\ref{pt3.3}) and (\ref{pt3.4}), we arrive at
\begin{align*}
&\big|\hat{K_1}(t,\xi)\big| \lesssim t e^{-|\xi|^{2(\sigma- \delta)}t},\,\, \big|\hat{K_0}(t,\xi)\big| \lesssim e^{-c |\xi|^{2(\sigma- \delta)}t} \text{ for small } |\xi|, \\
&\big|\hat{K_1}(t,\xi)\big| \lesssim t e^{-|\xi|^{2\delta}t},\,\, \big|\hat{K_0}(t,\xi) \big| \lesssim e^{-|\xi|^{2\delta}t} \text{ for large } |\xi|,
\end{align*}
where $c$ is a suitable positive constant. Hence, we may conclude all the desired statements.
\end{proof}

\noindent From Proposition \ref{md3.1.9} the following statement follows immediately.
\bmd \label{md3.1.10}
The following estimates hold in $\R^n$ for any $n \ge 1$:
\begin{align*}
&\big\|F^{-1}\big(|\xi|^a \hat{K_0}\big)(t,\cdot)\big\|_{L^\ity}\lesssim \begin{cases}
t^{-\frac{n+a}{2\delta}} \text{ for } t\in (0,1], &\\
t^{-\frac{n+a}{2(\sigma-\delta)}} \text{ for } t\in[1,\ity), &
\end{cases}\\
&\big\|F^{-1}\big(|\xi|^a \hat{K_1}\big)(t,\cdot)\big\|_{L^\ity}\lesssim \begin{cases}
t^{1-\frac{n+a}{2\delta}} \text{ for } t\in (0,1], &\\
t^{1-\frac{n+a}{2(\sigma-\delta)}} \text{ for } t\in[1,\ity), &
\end{cases}
\end{align*}
for any non-negative number $a$.
\emd

\subsection{$L^r$ estimates}
From the statements of Propositions \ref{md3.1.6} and \ref{md3.1.10}, by applying an interpolation argument we may conclude the following $L^r$ estimates.
\bmd \label{md3.1.12}
The following estimates hold in $\R^n$ for any $n \ge 1$:
\begin{align*}
\big\|F^{-1}\big(|\xi|^a \hat{K_0}\big)(t,\cdot)\big\|_{L^r}&\lesssim
\begin{cases}
t^{-(2+[\frac{n}{2}])(\frac{\sigma}{2\delta}-1)\frac{1}{r} -\frac{n}{2\delta}(1-\frac{1}{r})-\frac{a}{2\delta}} \text{ if } t\in (0,1], &\\
t^{-\frac{n}{2(\sigma-\delta)}(1-\frac{1}{r})-\frac{a}{2(\sigma-\delta)}} \text{ if } t\in[1,\ity),&
\end{cases} \\
\big\|F^{-1}\big(|\xi|^a \hat{K_1}\big)(t,\cdot)\big\|_{L^r}&\lesssim
\begin{cases}
t^{1-(1+[\frac{n}{2}])(\frac{\sigma}{2\delta}-1)\frac{1}{r} -\frac{n}{2\delta}(1-\frac{1}{r})-\frac{a}{2\delta}} \text{ if } t\in (0,1], &\\
t^{1-\frac{n}{2(\sigma-\delta)}(1-\frac{1}{r})-\frac{a}{2(\sigma-\delta)}} \text{ if } t\in[1,\ity), &
\end{cases}
\end{align*}
for all $r \in [1,\ity]$ and any non-negative number $a$.
\emd

\renewcommand{\proofname}{Proof of Theorem 1.}
\begin{proof}
In order to prove the first statement, we only apply Young's convolution inequality and use the statements in Proposition \ref{md3.1.12}. Taking account of some estimates related to the time derivative, we note that
$$ \partial_t \hat{K_0}= -|\xi|^{2\sigma}\hat{K_1}\text{ and } \partial_t \hat{K_1}= \hat{K_0}- |\xi|^{2\delta}\hat{K_1}. $$
Then, applying again Young's convolution inequality and Proposition \ref{md3.1.12}, we may conclude the second statement. Hence, the proof of Theorem 1 is completed.
\end{proof}

\subsection{$L^{q}- L^{q}$ linear estimates with additional $L^{m}$ regularity for the data}
\bmd \label{md3.7}
Let $\delta\in \big(0,\f{\sigma}{2}\big)$ in (\ref{pt6.3}), $q\in (1,\ity)$ and $m\in [1,q)$. Then the Sobolev solutions to (\ref{pt6.3}) satisfy the $(L^m \cap L^q)-L^q$ estimates
\begin{align*}
\big\||D|^a u(t,\cdot)\big\|_{L^q}& \lesssim
\begin{cases}
t^{-(2+[\frac{n}{2}])(\frac{\sigma}{2\delta}-1)}\|u_0\|_{L^m \cap H^{a,q}}+ t^{1-(1+[\frac{n}{2}])(\frac{\sigma}{2\delta}-1)-\frac{a}{2\delta}}\|u_1\|_{L^m \cap H^{[a-\sigma]^+,q}} \text{ if } t\in (0,1], &\\
(1+t)^{-\frac{n}{2(\sigma-\delta)}(1-\frac{1}{r})-\frac{a}{2(\sigma-\delta)}} \|u_0\|_{L^m \cap H^{a,q}}\\
\qquad \quad+ (1+t)^{1-\frac{n}{2(\sigma-\delta)}(1-\frac{1}{r})-\frac{a}{2(\sigma-\delta)}}\|u_1\|_{L^m \cap H^{[a-\sigma]^+,q}} \text{ if } t\in[1,\ity), &
\end{cases}\\
\big\||D|^a u_t(t,\cdot)\big\|_{L^q}& \lesssim
\begin{cases}
t^{-(2+[\frac{n}{2}])(\frac{\sigma}{2\delta}-1)}\big(\|u_0\|_{L^m \cap H^{a+\sigma,q}}+\|u_1\|_{L^m \cap H^{a,q}}\big) \text{ if } t\in (0,1], &\\
(1+t)^{-\frac{n}{2(\sigma-\delta)}(1-\frac{1}{r})-\frac{a+2\delta}{2(\sigma-\delta)}} \|u_0\|_{L^m \cap H^{a+\sigma,q}}\\
\qquad \quad+ (1+t)^{1-\frac{n}{2(\sigma-\delta)}(1-\frac{1}{r})-\frac{a+2\delta}{2(\sigma-\delta)}}\|u_1\|_{L^m \cap H^{a,q}} \text{ if } t\in[1,\ity), &
\end{cases}
\end{align*}
and the $L^q-L^q$ estimates
\begin{align*}
\big\||D|^a u(t,\cdot)\big\|_{L^q}& \lesssim
\begin{cases}
t^{-(2+[\frac{n}{2}])(\frac{\sigma}{2\delta}-1)} \|u_0\|_{H^{\sigma,q}}+ t^{1-(1+[\frac{n}{2}])(\frac{\sigma}{2\delta}-1)-\frac{a}{2\delta}}\|u_1\|_{H^{[a-\sigma]^+,q}} \text{ if } t\in (0,1], &\\
(1+t)^{-\frac{a}{2(\sigma-\delta)}} \|u_0\|_{H^{a,q}}+ (1+t)^{1-\frac{a}{2(\sigma-\delta)}} \|u_1\|_{H^{[a-\sigma]^+,q}} \text{ if } t\in[1,\ity), &
\end{cases}\\
\big\||D|^a u_t(t,\cdot)\big\|_{L^q}& \lesssim
\begin{cases}
t^{-(2+[\frac{n}{2}])(\frac{\sigma}{2\delta}-1)}\big(\|u_0\|_{H^{a+\sigma,q}}+\|u_1\|_{H^{a,q}}\big) \text{ if } t\in (0,1], &\\
(1+t)^{-\frac{a+2\delta}{2(\sigma-\delta)}} \|u_0\|_{H^{a+\sigma,q}}+ (1+t)^{1-\frac{a+2\delta}{2(\sigma-\delta)}}\|u_1\|_{H^{a,q}} \text{ if } t\in[1,\ity), &
\end{cases}
\end{align*}
where $1+ \frac{1}{q}= \frac{1}{r}+\frac{1}{m}$, for any non-negative number $a$ and for all $n \ge 1$.
\emd

\renewcommand{\proofname}{Proof.}
\begin{proof}
In order to obtain the $(L^m \cap L^q)- L^q$ estimates, we estimate the $L^q$ norm of the low-frequency part of the solutions by the $L^m$ norm of the data, whereas its high-frequency part is estimated by using the $L^q-L^q$ estimates. Thanks to Propositions \ref{md3.1.2} and \ref{md3.1.9}, we derive
\begin{align*}
\big\|F^{-1}\big(|\xi|^a \hat{K_0}\chi(|\xi|)\big)(t,\cdot)\big\|_{L^r}&\lesssim
\begin{cases}
1 \text{ if } t\in (0,1], &\\
t^{-\frac{n}{2(\sigma-\delta)}(1-\frac{1}{r})-\frac{a}{2(\sigma-\delta)}} \text{ if } t\in[1,\ity),&
\end{cases} \\
\big\|F^{-1}\big(|\xi|^a \hat{K_1}\chi(|\xi|)\big)(t,\cdot)\big\|_{L^r}&\lesssim
\begin{cases}
t \text{ if } t\in (0,1], &\\
t^{1-\frac{n}{2(\sigma-\delta)}(1-\frac{1}{r})-\frac{a}{2(\sigma-\delta)}} \text{ if } t\in[1,\ity), &
\end{cases}
\end{align*}
for all $r \in [1,\ity]$ and any non-negative number $a$. Therefore, applying Young's convolution inequality and using the suitable regularity of the data $u_0$ and $u_1$ depending on the order of $a$, we may conclude all the desired estimates for the solution and some its derivatives.
\end{proof}

\begin{nx}
\fontshape{n}
\selectfont
The singular behavior of the time-dependent coefficients for $t\longrightarrow +0$ in the above estimates brings some difficulties in the treatment of the semi-linear models (\ref{pt6.1}) and (\ref{pt6.2}). However, to avoid this difficulty in the proof of global (in time) existence results, we can compensate this singular behavior by assuming additional regularity of order $s_0$ for the data $u_0$ and $u_1$. For this reason, there appears a loss of regularity with respect to the initial data. We obtain the following corollary.
\end{nx}

\begin{hq} \label{hq6.6.1}
Let $\delta\in \big(0,\f{\sigma}{2}\big)$ in (\ref{pt6.3}), $q\in (1,\ity)$ and $m\in [1,q)$. Then the solution to (\ref{pt6.3}) satisfies the following $(L^m \cap L^q)-L^q$ estimates:
\begin{align*}
\big\||D|^a u(t,\cdot)\big\|_{L^q} &\lesssim (1+t)^{-\frac{n}{2(\sigma-\delta)}(1-\frac{1}{r})-\frac{a}{2(\sigma-\delta)}}\|u_0\|_{L^m \cap H^{a+s_0,q}}\\
&\qquad \quad+ (1+t)^{1-\frac{n}{2(\sigma-\delta)}(1-\frac{1}{r})-\frac{a}{2(\sigma-\delta)}}\|u_1\|_{L^m \cap H^{[a-\sigma+s_0]^+,q}},
\end{align*}
\begin{align*}
\big\||D|^{a}u_t(t,\cdot)\big\|_{L^q} &\lesssim (1+t)^{-\frac{n}{2(\sigma-\delta)}(1-\frac{1}{r})-\frac{a+2\delta}{2(\sigma-\delta)}}\|u_0\|_{L^m \cap H^{a+\sigma+s_0,q}}\\
&\qquad \quad+ (1+t)^{1-\frac{n}{2(\sigma-\delta)}(1-\frac{1}{r})-\frac{a+2\delta}{2(\sigma-\delta)}}\|u_1\|_{L^m \cap H^{a+s_0,q}},
\end{align*}
where $1+\frac{1}{q}=\frac{1}{r}+\frac{1}{m}$, for any $t>0$ and for all $n \ge 1$.
\end{hq}

%=================================================================================={Semi-linear estimates}
\section{Proof of the global (in time) existence results}
\label{Semi-linear estimates}

%=================================================================================={Philosophy of our approach}
\subsection{Philosophy of our approach}
In this section, we will apply the estimates for the solutions to (\ref{pt6.3}) from Proposition \ref{md3.7} to prove the global (in time) existence of small data solutions to the semi-linear models (\ref{pt6.1}) and (\ref{pt6.2}). By using the fundamental solutions $K_0$ and $K_1$ defined in Section $3$ we write the solutions to (\ref{pt6.3}) in the form
$$ u^{ln}(t,x)=K_0(t,x) \ast_{x} u_0(x)+ K_1(t,x) \ast_{x} u_1(x). $$
Since we are interested in dealing with semi-linear models with constant coefficients in the linear part, applying Duhamel's principle leads to the following formal implicit representation of the solutions to (\ref{pt6.1}) and (\ref{pt6.2}):
$$ u(t,x)= K_0(t,x) \ast_x u_0(x)+ K_1(t,x) \ast_x u_1(x) + \int_0^t K_1(t-\tau,x) \ast_x f(u,u_t) d\tau=: u^{ln}(t,x)+ u^{nl}(t,x), $$
where $f(u,u_t)=|u|^p$ or $|u_t|^p$. We choose the data space $\mathcal{A}$ and introduce the family $\{X(t)\}_{t>0}$ of solution spaces $X(t)$ with the norm
\begin{align*}
\|u\|_{X(t)}:= & \sup_{0\le \tau \le t} \Big( f_{1}(\tau)^{-1}\|u(\tau,\cdot)\|_{L^q}+ f_{2,s}(\tau)^{-1}\big\||D|^s u(\tau,\cdot)\big\|_{L^q}\\
&\qquad \quad + f_{3}(\tau)^{-1}\|u_t(\tau,\cdot)\|_{L^q}+ f_{4,s}(\tau)^{-1}\big\||D|^{s-\sigma} u_t(\tau,\cdot)\big\|_{L^q} \Big).
\end{align*}
Moreover, we introduce the family $\{X_0(t)\}_{t>0}$ of space $X_0(t):= C([0,t],H^{s,q})$ with the norm
$$\|u\|_{X_0(t)}:= \sup_{0\le \tau \le t} \Big( f_{1}(\tau)^{-1}\|u(\tau,\cdot)\|_{L^q}+ f_{2,s}(\tau)^{-1}\big\||D|^s u(\tau,\cdot)\big\|_{L^q} \Big), $$
where
\begin{align}
&f_{1}(\tau)= (1+\tau)^{1-\frac{n}{2(\sigma-\delta)}(1-\frac{1}{r})},\,\, f_{2,s}(\tau)=(1+\tau)^{1-\frac{n}{2(\sigma-\delta)}(1-\frac{1}{r})-\frac{s}{2(\sigma-\delta)}}, \label{pt4.1}\\
&f_{3}(\tau)=(1+\tau)^{1-\frac{n}{2(\sigma-\delta)}(1-\frac{1}{r})-\frac{\delta}{\sigma-\delta}},\,\, f_{4,s}(\tau)=(1+\tau)^{1-\frac{n}{2(\sigma-\delta)}(1-\frac{1}{r})-\frac{s-\sigma+2\delta}{2(\sigma-\delta)}}. \label{pt4.2}
\end{align}
We define for all $t>0$ the operator $N: \quad u \in X(t) \longrightarrow Nu \in X(t)$ by the formula
$$Nu(t,x)= K_0(t,x) \ast_x u_0(x)+ K_1(t,x) \ast_x u_1(x)+ \int_0^t K_1(t-\tau,x) \ast_x f(u,u_t) d\tau. $$
We will prove that the operator $N$ satisfies the following two inequalities:
\begin{align}
\|Nu\|_{X(t)}& \lesssim \|(u_0,u_1)\|_{\mathcal{A}}+ \|u\|^p_{X_0(t)}, \label{pt4.3}\\
\|Nu-Nv\|_{X(t)}& \lesssim \|u-v\|_{X_0(t)} \big(\|u\|^{p-1}_{X_0(t)}+ \|v\|^{p-1}_{X_0(t)}\big). \label{pt4.4}
\end{align}
Then, we apply Banach's fixed point theorem to gain local (in time) and global (in time) existence results as well.
\begin{nx}
\fontshape{n}
\selectfont
We can see that by plugging $a=s$ in the statements from Corollary \ref{hq6.6.1}, from the definition of the norm in $X(t)$ we conclude
\begin{equation}
\big\|u^{ln} \big\|_{X(t)} \lesssim \|(u_0,u_1)\|_{\mathcal{A}^{s+s_0}_{m,q}}, \text{ for }s \ge 0. \label{pt4.5}
\end{equation}
For this reason, to complete the proof of (\ref{pt4.3}) we need to show the following inequality:
\begin{equation}
\big\|u^{nl} \big\|_{X(t)} \lesssim \|u\|^p_{X_0(t)}. \label{pt4.31}
\end{equation}
\end{nx}
Now we are ready to prove our theorems from Section $2.2$.

%=================================================================================={No loss of decay}
\subsection{No loss of decay}

%==================================================================================={Theorem 2-A}
\renewcommand{\proofname}{Proof of Theorem 2-A: $s=\sigma$.}
\begin{proof}
We introduce the data space $\mathcal{A}:= \mathcal{A}^{\sigma+s_0}_{m,q}$ and the solution space
$$X(t):= C([0,t],H^{\sigma,q}) \cap C^1([0,t],L^q), $$
where the weight $f_{4,s}(\tau) \equiv 0$. First, let us prove the inequality (\ref{pt4.31}). In order to control $u^{nl}$, we use two different strategies for $\tau \in \big[0,[t-1]^+ \big]$ and $\tau \in \big[[t-1]^+,t \big]$. In particular, we use the $(L^m \cap L^q)- L^q$ estimates if $\tau \in \big[0,[t-1]^+ \big]$ and the $L^q-L^q$ estimates if $\tau \in \big[[t-1]^+,t \big]$ from Proposition \ref{md3.7}. Therefore, we derive for $j,k=0,1$ and $(j,k) \neq (1,1)$ the following estimates:
\begin{align*}
\big\|\partial_t^j |D|^{k\sigma} u^{nl}(t,\cdot)\big\|_{L^q} &\lesssim \int_0^{[t-1]^+}(1+t-\tau)^{1-\frac{n}{2(\sigma-\delta)}(1-\frac{1}{r})-\frac{k\sigma+2j\delta}{2(\sigma-\delta)}}\big\||u(\tau,\cdot)|^p\big\|_{L^m \cap L^q}d\tau\\
&\qquad + \int_{[t-1]^+}^t (t-\tau)^{1-(1+[\frac{n}{2}])(\frac{\sigma}{2\delta}-1)-(k+j)\frac{\sigma}{2\delta}}\big\||u(\tau,\cdot)|^p\big\|_{L^q}d\tau.
\end{align*}
Hence, it is necessary to require the estimates for $|u(\tau,x)|^p$ in $L^m \cap L^q$ and $L^q$ as follows:
$$\big\||u(\tau,\cdot)|^p\big\|_{L^m \cap L^q} \lesssim \|u(\tau,\cdot)\|^p_{L^{mp}}+ \|u(\tau,\cdot)\|^p_{L^{qp}},\, \text{ and }\big\||u(\tau,\cdot)|^p\big\|_{L^q}= \|u(\tau,\cdot)\|^p_{L^{qp}}.$$
Applying the fractional Gagliardo-Nirenberg inequality from Proposition \ref{fractionalGagliardoNirenberg} we can conclude
\begin{align*}
\big\||u(\tau,\cdot)|^p\big\|_{L^m \cap L^q} &\lesssim (1+\tau)^{p-\frac{n}{2m(\sigma-\delta)}{(p-1)}}\|u\|^p_{X_0(\tau)},\\
\big\||u(\tau,\cdot)|^p\big\|_{L^q} &\lesssim (1+\tau)^{p-\frac{np}{2(\sigma-\delta)}(\frac{1}{m}-\frac{1}{qp})}\|u\|^p_{X_0(\tau)},
\end{align*}
provided that (\ref{GN2A}) is satisfied. From both estimates we may conclude
\begin{align*}
\big\|\partial_t^j |D|^{k\sigma}u^{nl}(t,\cdot)\big\|_{L^q} &\lesssim \|u\|^p_{X_0(t)} \int_0^{[t-1]^+}(1+t-\tau)^{1-\frac{n}{2(\sigma-\delta)}(1-\frac{1}{r})-\frac{k\sigma+2j\delta}{2(\sigma-\delta)}}(1+\tau)^{p-\frac{n}{2m(\sigma-\delta)}{(p-1)}} d\tau\\
&\qquad + \|u\|^p_{X_0(t)} \int_{[t-1]^+}^t (t-\tau)^{1-(1+[\frac{n}{2}])(\frac{\sigma}{2\delta}-1)-(k+j)\frac{\sigma}{2\delta}}(1+\tau)^{p-\frac{np}{2(\sigma-\delta)}(\frac{1}{m}-\frac{1}{qp})}d\tau.
\end{align*}
The key tool relies now in Lemma \ref{LemmaIntegral}. Because of the condition (\ref{exponent2A}), applying Lemma \ref{LemmaIntegral} by choosing
$\alpha= -1+\frac{n}{2(\sigma-\delta)}\big(1-\frac{1}{r}\big)+\frac{k\sigma+2j\delta}{2(\sigma-\delta)}$ and
$\beta= -p+\frac{n}{2m(\sigma-\delta)}{(p-1)}$ we get
$$\int_0^{[t-1]^+}(1+t-\tau)^{1-\frac{n}{2(\sigma-\delta)}(1-\frac{1}{r})-\frac{k\sigma+2j\delta}{2(\sigma-\delta)}}(1+\tau)^{p-\frac{n}{2m(\sigma-\delta)}{(p-1)}} d\tau \lesssim (1+t)^{1-\frac{n}{2(\sigma-\delta)}(1-\frac{1}{r})-\frac{k\sigma+2j\delta}{2(\sigma-\delta)}}. $$
Moreover, since the condition $\big[\frac{n}{2}\big]< n_0$ holds, it follows $$ 1-\Big(1+\Big[\frac{n}{2}\Big]\Big)\Big(\frac{\sigma}{2\delta}-1\Big)-(k+j)\frac{\sigma}{2\delta}>-1, $$
and
$$ \int_{[t-1]^+}^t (t-\tau)^{1-(1+[\frac{n}{2}])(\frac{\sigma}{2\delta}-1)-(k+j)\frac{\sigma}{2\delta}}d\tau \lesssim \int_0^1 r^{1-(1+[\frac{n}{2}])(\frac{\sigma}{2\delta}-1)-(k+j)\frac{\sigma}{2\delta}}dr \lesssim 1.$$
Therefore, we can estimate
$$ \int_{[t-1]^+}^t (t-\tau)^{1-(1+[\frac{n}{2}])(\frac{\sigma}{2\delta}-1)-(k+j)\frac{\sigma}{2\delta}}(1+\tau)^{p-\frac{np}{2(\sigma-\delta)}(\frac{1}{m}-\frac{1}{qp})}d\tau \lesssim (1+t)^{1-\frac{n}{2(\sigma-\delta)}(1-\frac{1}{r})-\frac{k\sigma+2j\delta}{2(\sigma-\delta)}}. $$
Finally, we conclude for $j,k=0,1$ and $(j,k) \neq (1,1)$ the following estimate:
$$\big\|\partial_t^j |D|^{k\sigma}u^{nl}(t,\cdot)\big\|_{L^q} \lesssim (1+t)^{1-\frac{n}{2(\sigma-\delta)}(1-\frac{1}{r})-\frac{k\sigma+2j\delta}{2(\sigma-\delta)}} \|u\|^p_{X_0(t)}. $$
From the definition of the norm in $X(t)$, we obtain immediately the inequality (\ref{pt4.31}). \medskip

\noindent Next, let us prove the estimate (\ref{pt4.4}). Using again the $(L^m \cap L^q)- L^q$ estimates if $\tau \in \big[0,[t-1]^+ \big]$ and the $L^q-L^q$ estimates if $\tau \in \big[[t-1]^+,t \big]$ from Proposition \ref{md3.7}, we derive for two functions $u$ and $v$ from $X(t)$ the following estimate:
\begin{align*}
\big\|\partial_t^j |D|^{k\sigma}\big(Nu(t,\cdot)-Nv(t,\cdot)\big)\big\|_{L^q} &\lesssim \int_0^{[t-1]^+}(1+t-\tau)^{1-\frac{n}{2(\sigma-\delta)}(1-\frac{1}{r})-\frac{k\sigma+2j\delta}{2(\sigma-\delta)}}\big\||u(\tau,\cdot)|^p-|v(\tau,\cdot)|^p\big\|_{L^m \cap L^q}d\tau\\
&\qquad + \int_{[t-1]^+}^t (t-\tau)^{1-(1+[\frac{n}{2}])(\frac{\sigma}{2\delta}-1)-(k+j)\frac{\sigma}{2\delta}}\big\||u(\tau,\cdot)|^p-|v(\tau,\cdot)|^p\big\|_{L^q}d\tau.
\end{align*}
By using H\"{o}lder's inequality, we get
\begin{align*}
\big\||u(\tau,\cdot)|^p-|v(\tau,\cdot)|^p\big\|_{L^q}& \lesssim \|u(\tau,\cdot)-v(\tau,\cdot)\|_{L^{qp}} \big(\|u(\tau,\cdot)\|^{p-1}_{L^{qp}}+\|v(\tau,\cdot)\|^{p-1}_{L^{qp}}\big),\\
\big\||u(\tau,\cdot)|^p-|v(\tau,\cdot)|^p\big\|_{L^m}& \lesssim \|u(\tau,\cdot)-v(\tau,\cdot)\|_{L^{mp}} \big(\|u(\tau,\cdot)\|^{p-1}_{L^{mp}}+\|v(\tau,\cdot)\|^{p-1}_{L^{mp}}\big).
\end{align*}
Analogously to the proof of (\ref{pt4.3}), applying the fractional Gagliardo-Nirenberg inequality from Proposition \ref{fractionalGagliardoNirenberg} to the terms
$$ \|u(\tau,\cdot)-v(\tau,\cdot)\|_{L^\eta }, \text{ }\|u(\tau,\cdot)\|_{L^\eta}, \text{ }\|v(\tau,\cdot)\|_{L^\eta} $$
with $\eta=qp$ or $\eta=mp$ we may conclude the inequality (\ref{pt4.4}). Summarizing, the proof of Theorem 2-A is completed.
\end{proof}

%==================================================================================Theorem 3-A
\renewcommand{\proofname}{Proof of Theorem 3-A: $0 < s <\sigma$.}
\begin{proof}
We introduce the data space $\mathcal{A}:= \mathcal{A}^{s+s_0}_{m,q}$ and the solution space
$$X(t):= C([0,t],H^{s,q}), $$
where the weights $f_{3}(\tau)=f_{4,s}(\tau)\equiv 0$. We can see that $X_0(t)$ and $X(t)$ coincide in (\ref{pt4.4}) and (\ref{pt4.31}). In order to prove these two inequalities, we use the $(L^m \cap L^q)- L^q$ estimates if $\tau \in \big[0,[t-1]^+ \big]$ and the $L^q-L^q$ estimates if $\tau \in \big[[t-1]^+,t \big]$ from Proposition \ref{md3.7}. Therefore, we obtain for $k=0,1$ the following estimates:
\begin{align*}
\big\||D|^{ks}u^{nl}(t,\cdot)\big\|_{L^q} &\lesssim \int_0^{[t-1]^+}(1+t-\tau)^{1-\frac{n}{2(\sigma-\delta)}(1-\frac{1}{r})-\frac{ks}{2(\sigma-\delta)}}\big\||u(\tau,\cdot)|^p\big\|_{L^m \cap L^q}d\tau\\
&\qquad + \int_{[t-1]^+}^t (t-\tau)^{1-(1+[\frac{n}{2}])(\frac{\sigma}{2\delta}-1)-\frac{ks}{2\delta}}\big\||u(\tau,\cdot)|^p\big\|_{L^q}d\tau,
\end{align*}
 and
\begin{align*}
\big\||D|^{ks}\big(Nu(t,\cdot)-Nv(t,\cdot)\big)\big\|_{L^q} &\lesssim \int_0^{[t-1]^+}(1+t-\tau)^{1-\frac{n}{2(\sigma-\delta)}(1-\frac{1}{r})-\frac{ks}{2(\sigma-\delta)}}\big\||u(\tau,\cdot)|^p-|v(\tau,\cdot)|^p\big\|_{L^m \cap L^q}d\tau\\
&\qquad + \int_{[t-1]^+}^t (t-\tau)^{1-(1+[\frac{n}{2}])(\frac{\sigma}{2\delta}-1)-\frac{ks}{2\delta}}\big\||u(\tau,\cdot)|^p-|v(\tau,\cdot)|^p\big\|_{L^q}d\tau.
\end{align*}
In an analogous way as we did in the proof of Theorem 2-A, we may conclude for $k=0,1$ the following estimates:
\begin{align*}
\big\||D|^{ks}u^{nl}(t,\cdot)\big\|_{L^q} &\lesssim (1+t)^{1-\frac{n}{2(\sigma-\delta)}(1-\frac{1}{r})-\frac{ks}{2(\sigma-\delta)}}\|u\|^p_{X(t)},\\
\big\||D|^{ks}\big(Nu(t,\cdot)-Nv(t,\cdot)\big)\big\|_{L^q} &\lesssim (1+t)^{1-\frac{n}{2(\sigma-\delta)}(1-\frac{1}{r})-\frac{ks}{2(\sigma-\delta)}}\|u-v\|_{X(t)}\big(\|u\|^{p-1}_{X(t)}+\|v\|^{p-1}_{X(t)}\big),
\end{align*}
provided that the conditions (\ref{exponent3A}) and (\ref{GN3A}) are satisfied. From the definition of the norm in $X(t)$ we obtain immediately the inequalities (\ref{pt4.31}) and (\ref{pt4.4}). Summarizing, the proof of Theorem 3-A is completed.
\end{proof}

%==================================================================================={Theorem 4-A}
\renewcommand{\proofname}{Proof of Theorem 4-A: $\sigma< s \le \sigma+ \frac{n}{q}$.}
\begin{proof}
We introduce the data space $\mathcal{A}:= \mathcal{A}^{s+s_0}_{m,q}$ and the solution space
$$X(t):= C([0,t],H^{s,q}) \cap C^1([0,t],H^{s-\sigma,q}). $$
\noindent First, let us prove the inequality (\ref{pt4.31}). We have to control the norms
$$\|u^{nl}(t,\cdot)\|_{L^q},\,\, \|u_t^{nl}(t,\cdot)\|_{L^q},\,\, \big\||D|^{s}u^{nl}(t,\cdot)\big\|_{L^q},\,\, \big\||D|^{s-\sigma}u_t^{nl}(t,\cdot)\big\|_{L^q}. $$
In the same way as we did in the proof of Theorem 2-A, we may conclude the following estimates:
\begin{align}
\|u^{nl}(t,\cdot)\|_{L^q} &\lesssim (1+t)^{1-\frac{n}{2(\sigma-\delta)}(1-\frac{1}{r})} \|u\|^p_{X_0(t)}, \label{t4A1}\\
\|\partial_t u^{nl}(t,\cdot)\|_{L^q} &\lesssim  (1+t)^{1-\frac{n}{2(\sigma-\delta)}(1-\frac{1}{r})-\frac{\delta}{\sigma-\delta}} \|u\|^p_{X_0(t)}, \label{t4A2}
\end{align}
provided that condition ({\ref{exponent4A}}) is fulfilled and
\begin{equation}
p\in \Big[\frac{q}{m},\ity \Big)  \text{ if } n\le qs, \text{ or }p \in \Big[\f{q}{m}, \f{n}{n-qs}\Big] \text{ if } n \in \Big(qs, \f{q^2 s}{q-m}\Big]. \label{t4A3}
\end{equation}
\noindent Now, let us turn to estimate the norm $\big\||D|^s u^{nl}(t,\cdot)\big\|_{L^q}$. We use the $(L^m \cap L^q)- L^q$ estimates if $\tau \in \big[0,[t-1]^+\big]$ and the $L^q- L^q$ estimates if $\tau \in \big[[t-1]^+,t \big]$ from Proposition \ref{md3.7} to derive
\begin{align*}
\big\||D|^s u^{nl}(t,\cdot)\big\|_{L^q} &\lesssim \int_0^{[t-1]^+} (1+t-\tau)^{1-\frac{n}{2(\sigma-\delta)}(1-\frac{1}{r})-\frac{s}{2(\sigma-\delta)}} \big\||u(\tau,\cdot)|^p\big\|_{L^m \cap L^q \cap \dot{H}^{s-\sigma, q}} d\tau\\
&\qquad \quad + \int_{[t-1]^+}^t (t-\tau)^{-(2+[\frac{n}{2}])(\frac{\sigma}{2\delta}-1)} \big\||u(\tau,\cdot)|^p\big\|_{ L^q \cap \dot{H}^{s-\sigma, q}} d\tau.
\end{align*}
The integrals with $\big\||u(\tau,\cdot)|^p\big\|_{L^m \cap L^q}$ and $\big\||u(\tau,\cdot)|^p\big\|_{L^ q}$ will be handled as we did to get (\ref{t4A1}). To estimate the integral with $\big\||u(\tau,\cdot)|^p\big\|_{\dot{H}^{s-\sigma, q}}$, we shall apply Proposition \ref{Propfractionalchainrulegeneral} for the fractional chain rule with $p> \lceil s-\sigma \rceil$ and Proposition \ref{fractionalGagliardoNirenberg} for the fractional Gagliardo-Nirenberg inequality. Therefore, we obtain
\begin{align*}
&\big\||u(\tau,\cdot)|^p\big\|_{\dot{H}^{s-\sigma, q}} \lesssim \|u(\tau,\cdot)\|^{p-1}_{L^{q_1}}\,\,\big\||D|^{s-\sigma}u(\tau,\cdot)\big\|_{L^{q_2}}\\
&\qquad \lesssim \|u(\tau,\cdot)\|^{(p-1)(1-\theta_{q_1})}_{L^q}\,\,\big\||D|^s u(\tau,\cdot)\big\|^{(p-1)\theta_{q_1}}_{L^q}\,\,\|u(\tau,\cdot)\|^{1-\theta_{q_2}}_{L^q}\,\,\big\||D|^s u(\tau,\cdot)\big\|^{\theta_{q_2}}_{L^q}\\
&\qquad \lesssim (1+\tau)^{p-\frac{np}{2(\sigma-\delta)}(\frac{1}{m}-\frac{1}{qp})-\frac{s-\sigma}{2(\sigma-\delta)}}\|u\|^p_{X_0(\tau)},
\end{align*}
where
$$\frac{p-1}{q_1}+\frac{1}{q_2}= \frac{1}{q},\,\, \theta_{q_1}= \frac{n}{s}\Big(\frac{1}{q}-\frac{1}{q_1}\Big) \in [0,1],\,\, \theta_{q_2}= \frac{n}{s}\Big(\frac{1}{q}-\frac{1}{q_2}+\frac{s-\sigma}{n}\Big) \in \Big[\frac{s-\sigma}{s},1\Big]. $$
These conditions imply the restriction
\begin{equation}
1<p\le 1+\frac{q\sigma}{n-qs} \text{ if } n>qs, \text{ or } p>1 \text{ if } n \le qs. \label{t4A4}
\end{equation}
As a result, we can conclude
\begin{equation}
\big\||D|^s u^{nl}(t,\cdot)\big\|_{L^q} \lesssim (1+t)^{1-\frac{n}{2(\sigma-\delta)}(1-\frac{1}{r})-\frac{s}{2(\sigma-\delta)}} \|u\|^p_{X_0(t)}. \label{t4A5}
\end{equation}
In the same way we also get
\begin{equation}
\big\||D|^{s-\sigma}u_t^{nl}(t,\cdot)\big\|_{L^q} \lesssim (1+t)^{1-\frac{n}{2(\sigma-\delta)}(1-\frac{1}{r})-\frac{s-\sigma+2\delta}{2(\sigma-\delta)}} \|u\|^p_{X_0(t)}. \label{t4A6}
\end{equation}
Summarizing, from (\ref{t4A1}) to (\ref{t4A2}), (\ref{t4A5}) to (\ref{t4A6}) and the definition of the norm in $X(t)$ we obtain immediately the inequality (\ref{pt4.31}).\medskip

\noindent Next, let us prove the inequality (\ref{pt4.4}). Following the proof of Theorem 2-A, the new difficulty which appears is to estimate the norm $\big\||u(\tau,\cdot)|^p-|v(\tau,\cdot)|^p\big\|_{\dot{H}^{s-\sigma,q}}$. By using the integral representation
$$ |u(\tau,x)|^p-|v(\tau,x)|^p=p\int_0^1 \big(u(\tau,x)-v(\tau,x)\big)G\big(\omega u(\tau,x)+(1-\omega)v(\tau,x)\big)d\omega, $$
where $G(u)=u|u|^{p-2}$, we obtain
$$\big\||u(\tau,\cdot)|^p-|v(\tau,\cdot)|^p\big\|_{\dot{H}^{s-\sigma,q}} \lesssim \int_0^1 \Big\||D|^{s-\sigma}\Big(\big(u(\tau,\cdot)-v(\tau,\cdot)\big)G\big(\omega u(\tau,\cdot)+(1-\omega)v(\tau,\cdot)\big)\Big)\Big\|_{L^q}d\omega. $$
Thanks to the fractional Leibniz formula from Proposition \ref{fractionalLeibniz}, we can proceed as follows:
\begin{align*}
\big\||u(\tau,\cdot)|^p-|v(\tau,\cdot)|^p\big\|_{\dot{H}^{s-\sigma,q}} &\lesssim \big\||D|^{s-\sigma}\big(u(\tau,\cdot)-v(\tau,\cdot)\big)\big\|_{L^{r_1}} \int_0^1 \big\|G\big(\omega u(\tau,\cdot)+(1-\omega)v(\tau,\cdot)\big)\big\|_{L^{r_2}}d\omega\\
&\qquad + \|u(\tau,\cdot)-v(\tau,\cdot)\|_{L^{r_3}}  \int_0^1 \big\||D|^{s-\sigma}G\big(\omega u(\tau,\cdot)+(1-\omega)v(\tau,\cdot)\big)\big\|_{L^{r_4}}d\omega\\
&\lesssim \big\||D|^{s-\sigma}\big(u(\tau,\cdot)-v(\tau,\cdot)\big)\big\|_{L^{r_1}} \Big(\|u(\tau,\cdot)\|^{p-1}_{L^{r_2 (p-1)}}+ \|v(\tau,\cdot)\|^{p-1}_{L^{r_2 (p-1)}}\Big)\\
&\qquad + \|u(\tau,\cdot)-v(\tau,\cdot)\|_{L^{r_3}}  \int_0^1 \big\||D|^{s-\sigma}G\big(\omega u(\tau,\cdot)+(1-\omega)v(\tau,\cdot)\big)\big\|_{L^{r_4}}d\omega,
\end{align*}
where
$$\frac{1}{r_1}+\frac{1}{r_2}= \frac{1}{r_3}+\frac{1}{r_4}= \frac{1}{q}.$$
Employing the fractional Gargliardo-Nirenberg inequality from Proposition \ref{fractionalGagliardoNirenberg} implies
\begin{align*}
\big\||D|^{s-\sigma}\big(u(\tau,\cdot)-v(\tau,\cdot)\big)\big\|_{L^{r_1}}&\lesssim \|u(\tau,\cdot)-v(\tau,\cdot)\|^{\theta_1}_{\dot{H}^{s,q}}\,\,\|u(\tau,\cdot)-v(\tau,\cdot)\|^{1-\theta_1}_{L^q},\\
\|u(\tau,\cdot)\|_{L^{r_2 (p-1)}}&\lesssim \|u(\tau,\cdot)\|^{\theta_2}_{\dot{H}^{s,q}}\,\,\|u(\tau,\cdot)\|^{1-\theta_2}_{L^q},\\
\|u(\tau,\cdot)-v(\tau,\cdot)\|_{L^{r_3}}&\lesssim \|u(\tau,\cdot)-v(\tau,\cdot)\|^{\theta_3}_{\dot{H}^{s,q}}\,\,\|u(\tau,\cdot)-v(\tau,\cdot)\|^{1-\theta_3}_{L^q},
\end{align*}
where
$$\theta_1= \frac{n}{s}\Big(\frac{1}{q}-\frac{1}{r_1}+\frac{s-\sigma}{n}\Big) \in \Big[\frac{s-\sigma}{s},1\Big],\,\, \theta_2= \frac{n}{s}\Big(\frac{1}{q}-\frac{1}{r_2(p-1)}\Big) \in [0,1],\,\, \theta_3= \frac{n}{s}\Big(\frac{1}{q}-\frac{1}{r_3}\Big) \in [0,1]. $$
Moreover, since $\omega \in [0,1]$ is a parameter, we may apply again the fractional chain rule from Proposition \ref{Propfractionalchainrulegeneral} with $p >1+ \lceil s-\sigma \rceil$ and the fractional Gagliardo-Nirenberg inequality from Proposition \ref{fractionalGagliardoNirenberg} to conclude
\begin{align*}
&\big\||D|^{s-\sigma}G\big(\omega u(\tau,\cdot)+(1-\omega)v(\tau,\cdot)\big)\big\|_{L^{r_4}}\\
&\qquad \lesssim \|\omega u(\tau,\cdot)+(1-\omega)v(\tau,\cdot)\|^{p-2}_{L^{r_5}}\,\, \big\||D|^{s-\sigma}\big(\omega u(\tau,\cdot)+(1-\omega)v(\tau,\cdot)\big)\big\|_{L^{r_6}}\\
&\qquad \lesssim \|\omega u(\tau,\cdot)+(1-\omega)v(\tau,\cdot)\|^{(p-2)\theta_5+\theta_6}_{\dot{H}^{s,q}}\,\, \|\omega u(\tau,\cdot)+(1-\omega)v(\tau,\cdot)\|^{(p-2)(1-\theta_5)+1-\theta_6}_{L^q},
\end{align*}
where
$$\frac{p-2}{r_5}+\frac{1}{r_6}= \frac{1}{r_4},\,\, \theta_5= \frac{n}{s}\Big(\frac{1}{q}-\frac{1}{r_5}\Big) \in [0,1],\,\, \theta_6= \frac{n}{s}\Big(\frac{1}{q}-\frac{1}{r_6}+\frac{s-\sigma}{n}\Big) \in \Big[\frac{s-\sigma}{s},1\Big]. $$
Hence, we derive
\begin{align*}
&\int_0^1 \big\||D|^{s-\sigma}G\big(\omega u(\tau,\cdot)+(1-\omega)v(\tau,\cdot)\big)\big\|_{L^{r_4}}d\omega\\
&\quad \lesssim \big(\|u(\tau,\cdot)\|_{\dot{H}^{s,q}}+\|v(\tau,\cdot)\|_{\dot{H}^{s,q}}\big)^{(p-2)\theta_5+\theta_6}\, \big(\|u(\tau,\cdot)\|_{L^q}+\|v(\tau,\cdot)\|_{L^q} \big)^{(p-2)(1-\theta_5)+1-\theta_6}.
\end{align*}
Therefore, we conclude
$$\big\||u(\tau,\cdot)|^p-|v(\tau,\cdot)|^p\big\|_{\dot{H}^{s-\sigma,q}} \lesssim (1+\tau)^{p-\frac{np}{2(\sigma-\delta)}(\frac{1}{m}-\frac{1}{qp})-\frac{s-\sigma}{2(\sigma-\delta)}}\|u-v\|_{X_0(\tau)}\big( \|u\|^{p-1}_{X_0(\tau)}+ \|v\|^{p-1}_{X_0(\tau)} \big), $$
where we note that
$$\theta_1+ (p-1)\theta_2= \theta_3+ (p-2)\theta_5+ \theta_6= \frac{n}{s}\Big(\frac{p-1}{q}+\frac{s-\sigma}{n}\Big). $$
Summarizing, we have shown the estimate
$$\big\||D|^s \big(Nu(t,\cdot)- Nv(t,\cdot)\big)\big\|_{L^q} \lesssim (1+t)^{1-\frac{n}{2(\sigma-\delta)}(1-\frac{1}{r})-\frac{s}{2(\sigma-\delta)}} \|u-v\|_{X_0(t)}\big( \|u\|^{p-1}_{X_0(t)}+ \|v\|^{p-1}_{X_0(t)} \big). $$
By the same treatment, we may derive
$$ \big\||D|^{s-\sigma} \partial_t \big(Nu(t,\cdot)- Nv(t,\cdot)\big)\big\|_{L^q} \lesssim (1+t)^{1-\frac{n}{2(\sigma-\delta)}(1-\frac{1}{r})-\frac{s-\sigma+2\delta}{2(\sigma-\delta)}} \|u-v\|_{X_0(t)}\big( \|u\|^{p-1}_{X_0(t)}+ \|v\|^{p-1}_{X_0(t)} \big). $$
From the definition of the norm in $X(t)$ we have completed the proof of (\ref{pt4.4}). Summarizing, Theorem 4-A is proved completely.
\end{proof}
\begin{nx}
\fontshape{n}
\selectfont
In this remark, we want to clarify the possibility to choose actually the parameters $q_1$, $q_2$, $r_1,\cdots, r_6$ and $\theta_1,\cdots,\theta_6$ as required in the proof to Theorem 4-A. First, let us see that we can choose $q_1$, $q_2$ such that $\frac{p-1}{q_1}+\frac{1}{q_2}= \frac{1}{q},\,  \theta_{q_1}=\frac{n}{s}\big(\frac{1}{q}-\frac{1}{q_1}\big)\in [0,1],\, \theta_{q_2}=\frac{n}{s}\big(\frac{1}{q}-\frac{1}{q_2}+\frac{s-\sigma}{n}\big)\in \big[\frac{s-\sigma}{s},1\big]$, thanks to the following conditions:
\begin{equation}
2\le p \le 1+\frac{q\sigma}{n-qs} \text{ if } n>qs, \text{ or }p \ge 2 \text{ if }n \le qs. \label{Conditionequation}
\end{equation}
Namely, we can describe the requirements on $\theta_{q_1}$ and $\theta_{q_2}$ in terms of conditions on $q_1$ and $q_2$ as $ \frac{1}{q_1} \in \big[\frac{1}{q}-\frac{s}{n},\frac{1}{q}\big]$ and $\frac{1}{q_2} \in \big[\frac{1}{q}-\frac{\sigma}{n},\frac{1}{q}\big]$. Combining the second condition on $q_2$ and the expression $\frac{1}{q_2}=\frac{1}{q}-\frac{p-1}{q_1}$ we can obtain the condition on $q_1$ as $\frac{1}{q_1} \le \frac{\sigma}{n(p-1)}$ since $p \ge 2$. Hence, in order to guarantee the existence of $q_1$ and $q_2$ it is sufficient to intersect the two condition intervals for $q_1$ to become non-empty. For this reason we have the following condition:
$$\frac{1}{q}-\frac{s}{n} \le \frac{\sigma}{n(p-1)}, $$
which implies immediately (\ref{Conditionequation}).\medskip

For the choice of $r_1$ and $r_2$ such that $\frac{1}{r_1}+\frac{1}{r_2}= \frac{1}{q},\, \theta_1=\frac{n}{s}\big(\frac{1}{q}-\frac{1}{r_1}+\frac{s-\sigma}{n}\big)\in \big[\frac{s-\sigma}{s},1\big],\, \theta_2=\frac{n}{s}\big(\frac{1}{q}-\frac{1}{r_2(p-1)}\big)\in [0,1]$, we can repeat exactly the above arguments to find (\ref{Conditionequation}) by $r_1$ in place of $q_2$ and $r_2$ in place of $\frac{q_1}{p-1}$.\medskip

Let us devote now to explain the existence of suitable parameters $r_3,\cdots, r_6$ and $\theta_3,\cdots, \theta_6$. In the first step, our goal is to clarify $r_3$ and $r_4$ such that $\frac{1}{r_3}+\frac{1}{r_4}= \frac{1}{q}$ and $\theta_3=\frac{n}{s}\big(\frac{1}{q}-\frac{1}{r_3}\big) \in [0,1]$. By re-writing $\frac{1}{r_4}=\frac{1}{q}-\frac{1}{r_3}$ we can express the condition on $\theta_3$ equivalent to the condition on $r_4$ as $\frac{1}{r_4} \in \big[0,\frac{s}{n}\big]$. Therefore, choosing $r_4$ in the above admissible range we can take $r_3$ to guarantee $\theta_3 \in [0,1]$. In the second step, taking account of the conditions  $\theta_5=\frac{n}{s}\big(\frac{1}{q}-\frac{1}{r_5}\big) \in [0,1]$ and $\theta_6=\frac{n}{s}\big(\frac{1}{q}-\frac{1}{r_6}+\frac{s-\sigma}{n}\big) \in \big[\frac{s-\sigma}{s},1\big]$,
we can re-write the conditions on $r_5$ and $r_6$ as $\frac{1}{r_5}\in \big[\frac{1}{q}-\frac{s}{n},\frac{1}{q}\big]$ and $\frac{1}{r_6}\in \big[\frac{1}{q}-\frac{\sigma}{n},\frac{1}{q}\big]$, respectively. Moreover, using the sum $\frac{1}{r_4}=\frac{p-2}{r_5}+\frac{1}{r_6}$ and the obtained condition $\frac{1}{r_4} \in [0,\frac{s}{n}]$ we can express the condition on $r_6$ in an equivalent way as the condition on $r_5$, in particular, $\frac{1}{r_6}\le \frac{s}{n}-(p-2)\frac{1}{r_5}$ since $p \ge 2$. Hence, in order to ensure that we get a non-empty range for the parameter $r_6$ we need to have the second condition for $r_5$ as $(p-2)\frac{1}{r_5}\le \frac{s+\sigma}{n}- \frac{1}{q}$. Finally, we have to check the conditions on $p$ and $n$ coming from the requirement on the non-empty admissible range for $r_5$, that is,
$$ (p-2) \Big(\frac{1}{q}-\frac{s}{n}\Big) \le \frac{s+\sigma}{n}- \frac{1}{q}, $$
which follows immediately again from (\ref{Conditionequation}). Summarizing, we have shown that (\ref{Conditionequation}) is sufficient to guarantee the possibility to choose suitable parameters $q_1$, $q_2$, $r_1,\cdots, r_6$ and $\theta_1,\cdots,\theta_6$ in the proof to Theorem 4-A.
\end{nx}

%===================================================================================Theorem 5-A
\renewcommand{\proofname}{Proof of Theorem 5-A: $s>\sigma+\frac{n}{q}$.}
\begin{proof}
We introduce both spaces for the data and the solutions as in Theorem 4-A. We can repeat exactly, on the one hand, the estimates of the terms $|u(\tau,\cdot)|^p$ and $|u(\tau,\cdot)|^p-|v(\tau,\cdot)|^p$ in $L^m$ and $L^q$ as we did in the proof to Theorem 4-A. On the other hand, let us devote to estimate the above terms in $\dot{H}^{s-\sigma, q}$ by using results on fractional powers and the fractional Sobolev embedding.\medskip

In the first step, let us begin with $\big\||u(\tau,\cdot)|^p\big\|_{\dot{H}^{s-\sigma, q}}$. We apply Corollary \ref{Corfractionalhomogeneous} for fractional powers with $s-\sigma \in \big(\frac{n}{q},p\big)$ and Corollary \ref{CorollaryEmbedding} with a suitable $s^* <\frac{n}{q}$. Therefore, we obtain
$$ \big\||u(\tau,\cdot)|^p\big\|_{\dot{H}^{s-\sigma, q}}\lesssim \|u(\tau,\cdot)\|_{\dot{H}^{s-\sigma, q}}\|u(\tau,\cdot)\|^{p-1}_{L^\ity} \lesssim \|u(\tau,\cdot)\|_{\dot{H}^{s-\sigma, q}}\big(\|u(\tau,\cdot)\|_{\dot{H}^{s^*, q}}+ \|u(\tau,\cdot)\|_{\dot{H}^{s-\sigma, q}}\big)^{p-1}. $$
Applying the fractional Gagliardo-Nirenberg inequality from Proposition \ref{fractionalGagliardoNirenberg} we have
\begin{align*}
\|u(\tau,\cdot)\|_{\dot{H}^{s-\sigma, q}} &\lesssim \|u(\tau,\cdot)\|^{1-\theta_1}_{L^q}\big\||D|^s u(\tau,\cdot)\big\|^{\theta_1}_{L^q} \lesssim (1+\tau)^{1-\frac{n}{2(\sigma-\delta)}(1-\frac{1}{r})-\frac{s-\sigma}{2(\sigma-\delta)}}\|u\|_{X_0(\tau)}, \\
\|u(\tau,\cdot)\|_{\dot{H}^{s^*, q}} &\lesssim \|u(\tau,\cdot)\|^{1-\theta_2}_{L^q}\big\||D|^s u(\tau,\cdot)\big\|^{\theta_2}_{L^q} \lesssim (1+\tau)^{1-\frac{n}{2(\sigma-\delta)}(1-\frac{1}{r})- \frac{s^*}{2(\sigma-\delta)}}\|u\|_{X_0(\tau)},
\end{align*}
where $\theta_1= 1- \frac{\sigma}{s}$ and $\theta_2= \frac{s^*}{s}$. Hence, we derive
\begin{align*}
&\big\||u(\tau,\cdot)|^p\big\|_{\dot{H}^{s-\sigma, q}}\lesssim (1+\tau)^{p(1-\frac{n}{2(\sigma-\delta)}(1-\frac{1}{r}))-\frac{s-\sigma}{2(\sigma-\delta)}- (p-1)\frac{s^*}{2(\sigma-\delta)}} \|u\|_{X_0(\tau)} \\
&\qquad \lesssim (1+\tau)^{p-\frac{n}{2m(\sigma-\delta)}(p-1)}\|u\|^p_{X_0(\tau)},
\end{align*}
if we choose $s^*= \frac{n}{q}-\e$ with a sufficiently small positive number $\e$. \\
Next, let us estimate the norm $\big\||u(\tau,\cdot)|^p-|v(\tau,\cdot)|^p\big\|_{\dot{H}^{s-\sigma,q}}$. Then, repeating the proof of Theorem 4-A and using the same treatment as in the first step, we get
$$\big\||u(\tau,\cdot)|^p-|v(\tau,\cdot)|^p\big\|_{\dot{H}^{s-\sigma,q}} \lesssim (1+\tau)^{p-\frac{n}{2m(\sigma-\delta)}(p-1)} \|u-v\|_{X_0(t)}\big( \|u\|^{p-1}_{X_0(t)}+ \|v\|^{p-1}_{X_0(t)} \big), $$
provided that $p>2$ and $p> 1+s-\sigma$. Summarizing, the proof of Theorem 5-A is completed.
\end{proof}

%=================================================================================={Theorem 6-A}
\renewcommand{\proofname}{Proof of Theorem 6-A: $s>\sigma+\frac{n}{q}$.}
\begin{proof}
We introduce the data space and the solution space as in Theorem 4-A. In the proof of this theorem, the space $X_0(t)$ is replaced by the space $X(t)$ in both inequalities (\ref{pt4.4}) and (\ref{pt4.31}). First, let us prove the inequality (\ref{pt4.31}). In order to control $u^{nl}$, we apply the $L^m \cap L^q- L^q$ estimates on the interval $\big[0,[t-1]^+ \big]$ and $L^q-L^q$ estimates on the interval $\big[[t-1]^+,t \big]$ from Proposition \ref{md3.7}. Hence, we get
\begin{align*}
\|u^{nl}(t,\cdot)\|_{L^q}&\lesssim \int_0^{[t-1]^+}(1+t-\tau)^{1-\frac{n}{2(\sigma-\delta)}(1-\frac{1}{r})}\big\||u_t(\tau,\cdot)|^p\big\|_{L^m \cap L^q}d\tau\\
&\quad + \int_{[t-1]^+}^t (1+t-\tau)^{1-(1+[\frac{n}{2}])(\frac{\sigma}{2\delta}-1)}\big\||u_t(\tau,\cdot)|^p\big\|_{L^q}d\tau.
\end{align*}
We have
$$\big\||u_t(\tau,\cdot)|^p\big\|_{L^m \cap L^q} \lesssim \|u_t(\tau,\cdot)\|^p_{L^{mp}}+ \|u_t(\tau,\cdot)\|^p_{L^{qp}},\, \text{ and }\big\||u_t(\tau,\cdot)|^p\big\|_{L^q}= \|u_t(\tau,\cdot)\|^p_{L^{qp}}.$$
Applying the fractional Gagliardo-Nirenberg inequality from Proposition \ref{fractionalGagliardoNirenberg} implies
\begin{align*}
\big\||u_t(\tau,\cdot)|^p\big\|_{L^m \cap L^q} &\lesssim (1+\tau)^{p(1-\frac{n}{2(\sigma-\delta)}(\frac{1}{m}-\frac{1}{mp})-\frac{\delta}{\sigma-\delta})}\|u\|^p_{X(\tau)},\\
\big\||u_t(\tau,\cdot)|^p\big\|_{L^q} &\lesssim (1+\tau)^{p(1-\frac{n}{2(\sigma-\delta)}(\frac{1}{m}-\frac{1}{qp})-\frac{\delta}{\sigma-\delta})}\|u\|^p_{X(\tau)},
\end{align*}
provided that $p \in \big[\frac{q}{m},\ity \big) \text{ since } s > \sigma+\frac{n}{q}$. Therefore, we derive
\begin{align*}
\|u^{nl}(t,\cdot)\|_{L^q} &\lesssim \|u\|^p_{X(t)} \int_0^{[t-1]^+}(1+t-\tau)^{1-\frac{n}{2(\sigma-\delta)}(1-\frac{1}{r})}
(1+\tau)^{p(1-\frac{n}{2(\sigma-\delta)}(\frac{1}{m}-\frac{1}{mp})-\frac{\delta}{\sigma-\delta})} d\tau\\
&\quad + \|u\|^p_{X(t)} \int_{[t-1]^+}^t (t-\tau)^{1-(1+[\frac{n}{2}])(\frac{\sigma}{2\delta}-1)}
(1+\tau)^{p(1-\frac{n}{2(\sigma-\delta)}(\frac{1}{m}-\frac{1}{qp})-\frac{\delta}{\sigma-\delta})}d\tau.
\end{align*}
Because of condition (\ref{exponent6A}), after applying Lemma \ref{LemmaIntegral} with
\[ \alpha= -1+\frac{n}{2(\sigma-\delta)}\Big(1-\frac{1}{r}\Big)\,\,\,\mbox{ and}\,\,\,
\beta= p\Big(-1+\frac{n}{2(\sigma-\delta)}\Big(\frac{1}{m}-\frac{1}{mp}\Big)+\frac{\delta}{\sigma-\delta}\Big),\]
we get
$$\int_0^{[t-1]^+}(1+t-\tau)^{1-\frac{n}{2(\sigma-\delta)}(1-\frac{1}{r})}
(1+\tau)^{p(1-\frac{n}{2(\sigma-\delta)}(\frac{1}{m}-\frac{1}{mp})-\frac{\delta}{\sigma-\delta})} d\tau \lesssim (1+t)^{1-\frac{n}{2(\sigma-\delta)}(1-\frac{1}{r})}. $$
Moreover, since the condition $\big[\frac{n}{2}\big]< n_0$ holds, it follows $1-\big(1+[\frac{n}{2}]\big)\big(\frac{\sigma}{2\delta}-1\big)>-1.$ Therefore, we can estimate
$$ \int_{[t-1]^+}^t (t-\tau)^{1-(1+[\frac{n}{2}])(\frac{\sigma}{2\delta}-1)}
(1+\tau)^{p(1-\frac{n}{2(\sigma-\delta)}(\frac{1}{m}-\frac{1}{qp})-\frac{\delta}{\sigma-\delta})}d\tau \lesssim (1+t)^{1-\frac{n}{2(\sigma-\delta)}(1-\frac{1}{r})}. $$
Hence, we arrive at the following estimate:
\begin{equation}
\|u^{nl}(t,\cdot)\|_{L^q} \lesssim (1+t)^{1-\frac{n}{2(\sigma-\delta)}(1-\frac{1}{r})} \|u\|^p_{X(t)}. \label{t6A1}
\end{equation}
In the same way, we also conclude
\begin{equation}
\|\partial_t u^{nl}(t,\cdot)\|_{L^q}\lesssim  (1+t)^{1-\frac{n}{2(\sigma-\delta)}(1-\frac{1}{r})-\frac{\delta}{\sigma-\delta}} \|u\|^p_{X(t)}. \label{t6A2}
\end{equation}
Now, let us devote to estimate the norm $\big\||D|^s u^{nl}(t,\cdot)\big\|_{L^q}$. We derive
\begin{align*}
\big\||D|^s u^{nl}(t,\cdot)\big\|_{L^q} &\lesssim \int_0^{[t-1]^+} (1+t-\tau)^{1-\frac{n}{2(\sigma-\delta)}(1-\frac{1}{r})-\frac{s}{2(\sigma-\delta)}} \big\||u_t(\tau,\cdot)|^p\big\|_{L^m \cap L^q \cap \dot{H}^{s-\sigma, q}} d\tau\\
&\qquad+ \int_{[t-1]^+}^t (t-\tau)^{-(2+[\frac{n}{2}])(\frac{\sigma}{2\delta}-1)} \big\||u_t(\tau,\cdot)|^p\big\|_{ L^q \cap \dot{H}^{s-\sigma, q}} d\tau.
\end{align*}
The integrals with $\big\||u_t(\tau,\cdot)|^p\big\|_{L^m \cap L^q}$ and $\big\||u_t(\tau,\cdot)|^p\big\|_{L^ q}$ will be handled as before to get (\ref{t6A1}). To estimate the integral with $\big\||u_t(\tau,\cdot)|^p\big\|_{\dot{H}^{s-\sigma, q}}$, we apply Corollary \ref{Corfractionalhomogeneous} for fractional powers with $s-\sigma \in \big(\frac{n}{q},p\big)$ and Corollary \ref{CorollaryEmbedding} with a suitable $s^* <\frac{n}{q}$. Therefore, we obtain
$$ \big\||u_t(\tau,\cdot)|^p\big\|_{\dot{H}^{s-\sigma, q}} \lesssim \|u_t(\tau,\cdot)\|_{\dot{H}^{s-\sigma, q}}\|u_t(\tau,\cdot)\|^{p-1}_{L^\ity} \lesssim \|u_t(\tau,\cdot)\|_{\dot{H}^{s-\sigma, q}}\big(\|u_t(\tau,\cdot)\|_{\dot{H}^{s^*, q}}+ \|u_t(\tau,\cdot)\|_{\dot{H}^{s-\sigma, q}}\big)^{p-1}. $$
After applying the fractional Gagliardo-Nirenberg inequality from Proposition \ref{fractionalGagliardoNirenberg} it follows
$$ \|u_t(\tau,\cdot)\|_{\dot{H}^{s^*, q}} \lesssim \|u_t(\tau,\cdot)\|^{1-\theta}_{L^q}\big\||D|^{s-\sigma} u_t(\tau,\cdot)\big\|^{\theta}_{L^q} \lesssim (1+\tau)^{1-\frac{n}{2(\sigma-\delta)}(1-\frac{1}{r})-\frac{\delta}{\sigma-\delta}- \frac{s^*}{2(\sigma-\delta)}}\|u\|_{X(\tau)}, $$
where $\theta= \frac{s^*}{s-\sigma}$. Hence, we derive
\begin{align*}
&\big\||u_t(\tau,\cdot)|^p\big\|_{\dot{H}^{s-\sigma, q}} \lesssim (1+\tau)^{p(1-\frac{n}{2(\sigma-\delta)}(1-\frac{1}{r})-\frac{\delta}{\sigma-\delta})-\frac{s-\sigma}{2(\sigma-\delta)}- (p-1)\frac{s^*}{2(\sigma-\delta)}}\|u\|^p_{X(\tau)}\\
&\qquad \lesssim (1+\tau)^{p(1-\frac{n}{2(\sigma-\delta)}(\frac{1}{m}-\frac{1}{mp})-\frac{\delta}{\sigma-\delta})}\|u\|^p_{X(\tau)},
\end{align*}
if we choose $s^*= \frac{n}{q}- \e$, where $\e$ is a sufficiently small positive. Analogous to the above arguments we may conclude
\begin{equation}
\big\||D|^s u^{nl}(t,\cdot)\big\|_{L^q} \lesssim (1+t)^{1-\frac{n}{2(\sigma-\delta)}(1-\frac{1}{r})-\frac{s}{2(\sigma-\delta)}} \|u\|^p_{X(t)}. \label{t6A3}
\end{equation}
In the same way, we also obtain
\begin{equation}
\big\||D|^{s-\sigma} u_t^{nl}(t,\cdot)\big\|_{L^q} \lesssim (1+t)^{1-\frac{n}{2(\sigma-\delta)}(1-\frac{1}{r})-\frac{s-\sigma+2\delta}{2(\sigma-\delta)}} \|u\|^p_{X(t)}. \label{t6A4}
\end{equation}
From (\ref{t6A1}) to (\ref{t6A4}) and the definition of the norm in $X(t)$ we obtain immediately the inequality (\ref{pt4.31}).\medskip

Next, let us prove the inequality (\ref{pt4.4}). The difficulty appearing is to cope with estimating the term $\big\||u_t(\tau,\cdot)|^p-|v_t(\tau,\cdot)|^p\big\|_{\dot{H}^{s-\sigma,q}}$. Then, repeating the proof of Theorem 4-A and using the same treatment as in the above first step, we get
$$ \big\||D|^s \big(Nu(t,\cdot)- Nv(t,\cdot)\big)\big\|_{L^q} \lesssim (1+t)^{1-\frac{n}{2(\sigma-\delta)}(1-\frac{1}{r})-\frac{s}{2(\sigma-\delta)}} \|u-v\|_{X(t)}\big( \|u\|^{p-1}_{X(t)}+ \|v\|^{p-1}_{X(t)} \big). $$
By the same treatment, we may conclude
$$ \big\||D|^{s-\sigma} \partial_t \big(Nu(t,\cdot)- Nv(t,\cdot)\big)\big\|_{L^q} \lesssim (1+t)^{1-\frac{n}{2(\sigma-\delta)}(1-\frac{1}{r})-\frac{s-\sigma+2\delta}{2(\sigma-\delta)}} \|u-v\|_{X(t)}\big( \|u\|^{p-1}_{X(t)}+ \|v\|^{p-1}_{X(t)} \big). $$
From the definition of the norm in $X(t)$ we have completed the proof of (\ref{pt4.4}). Summarizing, Theorem 6-A is proved completely.
\end{proof}

%=================================================================================={Loss of decay}
\subsection{Loss of decay}
In this section, we show how the restrictions to the admissible exponents $p$ appearing in all the theorems A can be relaxed. We will use some decay rates for solutions or some of their derivatives to the semi-linear models which are worse than those given for the solutions to the corresponding linear models with vanishing right-hand side to treat the semi-linear models (\ref{pt6.1}) and (\ref{pt6.2}), that is, we allow a loss of decay. This strategy comes into play to bring some advantage to weaken the restrictions to the admissible exponents $p$. In particular, we shall modify the weights in (\ref{pt4.1}) and (\ref{pt4.2}) to create a loss of decay in the following way:
$$ f_{1}(\tau):=f_{\e_1}(\tau)= (1+\tau)^{1-\frac{n}{2(\sigma-\delta)}(1-\frac{1}{r})+\e_1},\,\, f_{2,s}(\tau):=f_{\e_2,s}(\tau)= (1+\tau)^{1-\frac{n}{2(\sigma-\delta)}(1-\frac{1}{r})-\frac{s}{2(\sigma-\delta)}+\e_2}, $$
and $$ f_{3}(\tau):=f_{\e_3}(\tau)= (1+\tau)^{1-\frac{n}{2(\sigma-\delta)}(1-\frac{1}{r})-\frac{\delta}{\sigma-\delta}+\e_3},\,\, f_{4,s}(\tau):= f_{\e_4,s}(\tau)= (1+\tau)^{1-\frac{n}{2(\sigma-\delta)}(1-\frac{1}{r})-\frac{s-\sigma+2\delta}{2(\sigma-\delta)}+\e_4}, $$
for some positive constants $\e_j$ with $j=1,\cdots,4$. Here these constants stand for the loss of decay in comparison with the corresponding decay estimates for the solutions to (\ref{pt6.3}).\medskip

%=================================================================================={Theorem 2-B}
\renewcommand{\proofname}{Proof of Theorem 2-B: $s=\sigma$.}
\begin{proof}
We follow the proof of Theorem 2-A. Having in mind we fix the data space and the solution space as in Theorem 2-A, but we use the different weights where the weight $f_{\e_4,s}(\tau)\equiv 0$. In order to prove the inequality (\ref{pt4.3}), repeating the proof of Theorem 2-A we derive the following estimate:
\begin{align*}
&\big\|\partial_t^j |D|^{k\sigma} Nu(t,\cdot)\big\|_{L^q}\\
&\qquad \lesssim (1+t)^{1-\frac{n}{2(\sigma-\delta)}(1-\frac{1}{r})-\frac{k\sigma+2j\delta}{2(\sigma-\delta)}}
\|(u_0,u_1)\|_{\mathcal{A}^{\sigma+s_0}_{m,q}}\\
&\qquad \quad + \|u\|^p_{X_0(t)} \int_0^{[t-1]^+}(1+t-\tau)^{1-\frac{n}{2(\sigma-\delta)}(1-\frac{1}{r})-\frac{k\sigma+2j\delta}{2(\sigma-\delta)}}
(1+\tau)^{p-\frac{n}{2m(\sigma-\delta)}{(p-1)}+p(\e_2\theta_{mp}+\e_1(1-\theta_{mp}))} d\tau\\
&\qquad \quad + \|u\|^p_{X_0(t)} \int_{[t-1]^+}^t (t-\tau)^{1-(1+[\frac{n}{2}])(\frac{\sigma}{2\delta}-1)-(k+j)\frac{\sigma}{2\delta}}
(1+\tau)^{p-\frac{np}{2(\sigma-\delta)}(\frac{1}{m}-\frac{1}{qp})+p(\e_2\theta_{qp}+\e_1(1-\theta_{qp}))}d\tau.
\end{align*}
Now we fix the constant $\e_1:= \big(1-\frac{1}{p}\big)\big(-1+\frac{n}{2(\sigma-\delta)}\big(1-\frac{1}{r}\big)\big)$. Due to $n> n_1$, it follows $-1+\frac{n}{2(\sigma-\delta)}\big(1-\frac{1}{r}\big)>1$ and $\e_1$ is positive. Next we choose $\e_2=\frac{\sigma}{2(\sigma-\delta)}+\e_1$ and $\e_3=\frac{\delta}{\sigma-\delta}$. Then, we have
\begin{align*}
\big\|\partial_t^j |D|^{k\sigma} Nu(t,\cdot)\big\|_{L^q} &\lesssim (1+t)^{1-\frac{n}{2(\sigma-\delta)}(1-\frac{1}{r})-\frac{k\sigma+2j\delta}{2(\sigma-\delta)}}\|(u_0,u_1)\|_{\mathcal{A}^{\sigma+s_0}_{m,q}}\\
&\quad+ \|u\|^p_{X_0(t)} \int_0^{[t-1]^+}(1+t-\tau)^{1-\frac{n}{2(\sigma-\delta)}(1-\frac{1}{r})-\frac{k\sigma+2j\delta}{2(\sigma-\delta)}}(1+\tau)^{1-\frac{n}{2(\sigma-\delta)}(1-\frac{1}{r})} d\tau\\
&\quad + \|u\|^p_{X_0(t)} \int_{[t-1]^+}^t (t-\tau)^{1-(1+[\frac{n}{2}])(\frac{\sigma}{2\delta}-1)-(k+j)\frac{\sigma}{2\delta}}(1+\tau)^{1-\frac{n}{2(\sigma-\delta)}(1-\frac{1}{r})}d\tau.
\end{align*}
Applying Lemma \ref{LemmaIntegral} by choosing
$\alpha= -1+\frac{n}{2(\sigma-\delta)}\big(1-\frac{1}{r}\big)+\frac{k\sigma+2j\delta}{2(\sigma-\delta)}$
and $\beta= -1+\frac{n}{2(\sigma-\delta)}\big(1-\frac{1}{r}\big)$, we get
$$ \int_0^{[t-1]^+}(1+t-\tau)^{1-\frac{n}{2(\sigma-\delta)}(1-\frac{1}{r})-\frac{k\sigma+2j\delta}{2(\sigma-\delta)}}(1+\tau)^{1-\frac{n}{2(\sigma-\delta)}(1-\frac{1}{r})} d\tau \lesssim (1+t)^{1-\frac{n}{2(\sigma-\delta)}(1-\frac{1}{r})}. $$
Following the same arguments we used in the proof to Theorem 2-A, the condition $\big[\frac{n}{2}\big]< n_0$ implies
$$ \int_{[t-1]^+}^t (t-\tau)^{1-(1+[\frac{n}{2}])(\frac{\sigma}{2\delta}-1)-(k+j)\frac{\sigma}{2\delta}}(1+\tau)^{1-\frac{n}{2(\sigma-\delta)}(1-\frac{1}{r})}d\tau \lesssim (1+t)^{1-\frac{n}{2(\sigma-\delta)}(1-\frac{1}{r})}. $$
Finally, we conclude the following estimates:
\begin{align*}
f_{\e_1}(\tau)^{-1}\|Nu(\tau,\cdot)\|_{L^q}& \lesssim (1+t)^{-\e_1} \big( \|(u_0,u_1)\|_{\mathcal{A}}+ \|u\|^p_{X_0(t)}\big), \\
f_{\e_2,\sigma}(\tau)^{-1}\big\||D|^\sigma Nu(\tau,\cdot)\big\|_{L^q}& \lesssim (1+t)^{-\e_2} \|(u_0,u_1)\|_{\mathcal{A}}+ (1+t)^{-\e_1} \|u\|^p_{X_0(t)}, \\
f_{\e_3}(\tau)^{-1}\|\partial_t Nu(\tau,\cdot)\|_{L^q}& \lesssim (1+t)^{-\frac{\delta}{\sigma-\delta}} \|(u_0,u_1)\|_{\mathcal{A}}+ \|u\|^p_{X_0(t)}.
\end{align*}
From the definition of the norm in $X(t)$, we obtain immediately the inequality (\ref{pt4.3}). An analogous treatment as we did in the proof of Theorem 2-A and the above arguments give the following estimates:
\begin{align*}
f_{\e_1}(\tau)^{-1}\|Nu(\tau,\cdot)-Nv(\tau,\cdot)\|_{L^q}& \lesssim (1+t)^{-\e_1} \|u-v\|_{X_0(t)}\big(\|u\|^{p-1}_{X_0(t)}+\|v\|^{p-1}_{X_0(t)}\big), \\
f_{\e_2,\sigma}(\tau)^{-1}\big\||D|^\sigma \big(Nu(\tau,\cdot)-Nv(\tau,\cdot)\big)\big\|_{L^q}& \lesssim (1+t)^{-\e_1} \|u-v\|_{X_0(t)}\big(\|u\|^{p-1}_{X_0(t)}+\|v\|^{p-1}_{X_0(t)}\big), \\
f_{\e_3}(\tau)^{-1}\big\|\partial_t \big(Nu(\tau,\cdot)-Nv(\tau,\cdot)\big)\big\|_{L^q}& \lesssim \|u-v\|_{X_0(t)}\big(\|u\|^{p-1}_{X_0(t)}+\|v\|^{p-1}_{X_0(t)}\big).
\end{align*}
From the definition of the norm in $X(t)$, we obtain immediately the inequality (\ref{pt4.4}). Summarizing, the proof of Theorem 2-B is completed.
\end{proof}

%=================================================================================={Theorem 3-B}
\renewcommand{\proofname}{Proof of Theorem 3-B: $0 < s <\sigma$.}
\begin{proof}
We follow the proofs of Theorems 3-A and 2-B. Having in mind we fix the data space and the solution space as in Theorem 3-A, but we use the different weights where the weights $f_{\e_3}(\tau)= f_{\e_4,s}(\tau)\equiv 0$. Now we fix the constant $\e_1:= \big(1-\frac{1}{p}\big)\big(-1+\frac{n}{2(\sigma-\delta)}\big(1-\frac{1}{r}\big)\big)$. Next we choose $\e_2=\frac{s}{2(\sigma-\delta)}+\e_1$. Then, following the proofs of Theorems 3-A and 2-B we can prove Theorem 3-B.
\end{proof}

%=================================================================================={Theorem 4-B}
\renewcommand{\proofname}{Proof of Theorem 4-B: $\sigma< s \le \sigma+ \frac{n}{q}$.}
\begin{proof}
We follow the proofs of Theorems 4-A and 2-B. Having in mind we fix the data space and the solution space as in Theorem 4-A, but we use the different weights. Now we fix the constant $\e_1:= \big(1-\frac{1}{p}\big)\big(-1+\frac{n}{2(\sigma-\delta)}\big(1-\frac{1}{r}\big)\big)$. Next we choose $\e_2=\frac{s}{2(\sigma-\delta)}+\e_1$, $\e_3=\frac{\delta}{\sigma-\delta}$ and $\e_4=\frac{s-\sigma+2\delta}{2(\sigma-\delta)}$. Then, following the proofs of Theorems 4-A and 2-B we can prove Theorem 4-B.
\end{proof}

%=================================================================================={Theorem 5-B}
\renewcommand{\proofname}{Proof of Theorem 5-B: $s>\sigma+\frac{n}{q}$.}
\begin{proof}
We follow the proofs of Theorems 5-A and 2-B. Having in mind we fix both spaces of data and solution as in Theorem 5-A, but we use the different weights. Now we fix the constant $\e_1$, and choose $\e_j$ with $j=2,\cdots,4$ as in the proof of Theorem 4-B. Then, following the proofs of Theorems 5-A and 2-B we can prove Theorem 5-B.
\end{proof}

%=================================================================================={Theorem 6-B}
\renewcommand{\proofname}{Proof of Theorem 6-B: $s>\sigma+\frac{n}{q}$.}
\begin{proof}
We follow the proofs of Theorems 6-A and 2-B. Having in mind we fix the data space and the solution space as in Theorem 6-A, but we use the different weights. Now we fix the constant $\e:= \big(1-\frac{1}{p}\big)\big(-1+\frac{n}{2(\sigma-\delta)}\big(1-\frac{1}{r}\big)\big)$. Next we choose $\e_1=0$, $\e_2=\frac{s}{2(\sigma-\delta)}$, $\e_3=\frac{\delta}{\sigma-\delta}+\e$ and $\e_4=\frac{s-\sigma+2\delta}{2(\sigma-\delta)}+\e$. Then, following the proofs of Theorems 6-A and 2-B we can prove Theorem 6-B.
\end{proof}

%=================================================================================={Concluding remarks and open problems}
\section{Concluding remarks and open problems}

\begin{nx}{(Semi-linear visco-elastic damped $\sigma$-evolution models)}
\fontshape{n}
\selectfont
In this paper, we presented $(L^{m}\cap L^{q})- L^{q}$ and $L^{q}- L^{q}$ estimates for solutions and its derivatives to the model (\ref{pt6.3}) to prove the global (in time) existence of small data Sobolev solutions to the semi-linear models (\ref{pt6.1}) and (\ref{pt6.2}) with $\delta \in (0,\frac{\sigma}{2})$. It can be expected that the same approach could be applied to obtain $L^{1}$ estimates for oscillating integrals, simultaneously, $(L^{m}\cap L^{q})- L^{q}$ and $L^{q}- L^{q}$ estimates for solutions and its derivatives to the model (\ref{pt6.3}) with $\delta \in (\frac{\sigma}{2},\sigma)$ and the so-called visco-elastic type damped case $\delta=\sigma$ as well. We remark that the properties of the solutions to (\ref{pt6.1}) and (\ref{pt6.2}) change completely from $(0,\frac{\sigma}{2})$ to $(\frac{\sigma}{2},\sigma]$. In particular, we want to distinguish between the ``parabolic like models" ($\delta \in [0,\frac{\sigma}{2})$) and ``$\sigma$-evolution like models" ($\delta \in (\frac{\sigma}{2},\sigma]$) according to expected decay estimates. Hence, the following semi-linear models are of our interest:
\begin{equation}
u_{tt}+ (-\Delta)^\sigma u+ \mu (-\Delta)^\delta u_t = f(u,u_t),\, u(0,x)= u_0(x),\, u_t(0,x)=u_1(x) \label{pt0.1}
\end{equation}
with $\sigma \ge 1$, $\mu>0$ and $\delta \in (\frac{\sigma}{2},\sigma]$. In a forthcoming paper, we will study the global (in time) existence of small data Sobolev solutions from suitable spaces on the base of $L^q$ by assuming additional $L^{m}$ regularity for the initial data to (\ref{pt0.1}).
\end{nx}

\begin{nx}{(Time-dependent coefficients)}
\fontshape{n}
\selectfont
A next challenge is to obtain $(L^{m}\cap L^{q})- L^{q}$ and $L^{q}- L^{q}$ estimates for solutions to structurally damped $\sigma$-evolution models with time-dependent coefficients. These estimates are key tools to prove global (in time) existence results to semi-linear models. Hence, it is interesting to study the following Cauchy problem:
\begin{equation}
u_{tt}+ (-\Delta)^\sigma u+ b(t) (-\Delta)^\delta u_t = 0,\, u(0,x)= u_0(x),\, u_t(0,x)=u_1(x) \label{pt0.2}
\end{equation}
with $\sigma \ge 1$ and $\delta \in (0,\sigma)$. Here the coefficient $b=b(t)$ in (\ref{pt0.2}) should satisfy some ``effectiveness assumptions'' as in \cite{Kainane}.
\end{nx}

\begin{nx}{(Gevrey smoothing)}
\fontshape{n}
\selectfont
We are interested to understand to which Gevrey space the solutions to (\ref{pt6.3}) belong to. For this reason, we will use our estimates with $L^2$ norms and assume for the Cauchy data $(u_0,u_1) \in H^\sigma \times L^2$. The study of regularity properties for the solutions allows to restrict our considerations to large frequencies in the extended phase space. First, we recall the following definition of the Gevrey-Sobolev space regularity (see, for instance, \cite{ChenRodino,Kainane}).
%\bdn
%A given function $u: \R^n\longrightarrow \R$ belongs to the Gevrey space $\Gamma^{a,s}$ if and only if there exist positive constants $a$ and $s$ such that
%$$\exp \big(a \big<\xi\big>^{\frac{1}{s}}\big)F_{x\rightarrow \xi}(u)(\xi) \in L^2. $$
%Moreover, by $\Gamma^{s}$ we denote the inductive limit of all spaces $\Gamma^{a,s}$, that is, $\Gamma^{s}:= \bigcup_{a>0} \Gamma^{a,s}$.
%\edn

\bdn
A given function $u: \R^n\longrightarrow \R$ belongs to the Gevrey-Sobolev space $\Gamma^{a,s,\rho}$ if and only if there exist positive constants $a$, $s$ and a constant $\rho$ such that
$$\exp \big(a \big<\xi\big>^{\frac{1}{s}}\big) \big<\xi\big>^\rho F_{x\rightarrow \xi}(u)(\xi) \in L^2. $$
Moreover, by $\Gamma^{s,\rho}$ we denote the inductive limit of all spaces $\Gamma^{a,s,\rho}$, that is, $\Gamma^{s,\rho}:= \bigcup_{a>0} \Gamma^{a,s,\rho}$.
\edn
Then we may conclude the following statement.\medskip

\noindent\textbf{Theorem 7.}\textit{ Let us consider the Cauchy problem (\ref{pt6.3}) with $\delta\in (0,\frac{\sigma}{2})$. The data are supposed to belong to the energy space, that is, $(u_0,u_1) \in H^\sigma \times L^2$. Then, there is a smoothing effect in the sense, that the solution and its derivatives belong to the Gevrey-Sobolev space and the Gevrey space, respectively, that is,
\[ u(t,\cdot) \in \Gamma^{\frac{1}{2\delta},\sigma} \text{ and } |D|^\sigma u(t,\cdot),\,\, u_t(t,\cdot) \in \Gamma^{\frac{1}{2\delta},0} \text{ for all } t>0. \]}
\renewcommand{\proofname}{Proof.}
\begin{proof}
Using the asymptotic behavior of the characteristic roots for large $|\xi|$ in (\ref{pt3.3}) and (\ref{pt3.4}) we find
$$ |\hat{K_0}|\lesssim e^{-c|\xi|^{2\delta}t},\,\, |\hat{K_1}|\lesssim |\xi|^{-\sigma}e^{-c|\xi|^{2\delta}t} \text{ and } |\partial_t \hat{K_0}|\lesssim |\xi|^{\sigma}e^{-c|\xi|^{2\delta}t},\,\, |\partial_t \hat{K_1}|\lesssim e^{-c|\xi|^{2\delta}t}, $$
for some positive constant $c$. Hence, by the representation of the solutions (\ref{pt3.2}) we can easily see that the following relations hold for high frequencies:
\begin{align*}
\int_{\R^n}\exp \big(2c|\xi|^{2\delta}t \big) |\xi|^{2\sigma}(1-\chi(\xi))|v(t,\xi)|^2 d\xi &\lesssim \int_{\R^n}|\xi|^{2\sigma} |v_0(\xi)|^2 d\xi + \int_{\R^n}|v_1(\xi)|^2 d\xi, \\
\int_{\R^n}\exp \big(2c|\xi|^{2\delta}t \big) (1-\chi(\xi))|v_t(t,\xi)|^2 d\xi &\lesssim \int_{\R^n}|\xi|^{2\sigma} |v_0(\xi)|^2 d\xi + \int_{\R^n}|v_1(\xi)|^2 d\xi.
\end{align*}
Therefore, we may conclude immediately all the desired statements.
\end{proof}

\end{nx}





\noindent \textbf{Acknowledgments}\medskip

\noindent The PhD study of MSc. T.A. Dao is supported by Vietnamese Government's Scholarship.\medskip

%=================================================================================={Appendix}
\noindent\textbf{Appendix A}\medskip

%==================================================================
\noindent \textit{A.1. Fractional Gagliardo-Nirenberg inequality}

\begin{md} \label{fractionalGagliardoNirenberg}
Let $1<p,p_0,p_1<\infty$, $\sigma >0$ and $s\in [0,\sigma)$. Then, it holds the following fractional Gagliardo-Nirenberg inequality for all $u\in L^{p_0} \cap \dot{H}^\sigma_{p_1}$:
$$ \|u\|_{\dot{H}^{s}_p}\lesssim \|u\|_{L^{p_0}}^{1-\theta}\|u\|_{\dot{H}^{\sigma}_{p_1}}^\theta, $$
where $\theta=\theta_{s,\sigma}(p,p_0,p_1)=\frac{\frac{1}{p_0}-\frac{1}{p}+\frac{s}{n}}{\frac{1}{p_0}-\frac{1}{p_1}+\frac{\sigma}{n}}$ and $\frac{s}{\sigma}\leq \theta\leq 1$ .
\end{md}
For the proof one can see \cite{Ozawa}.
\medskip

%==================================================================
\noindent \textit{A.2. Fractional Leibniz rule}

\begin{md} \label{fractionalLeibniz}
Let us assume $s>0$ and $1\leq r \leq \infty, 1<p_1,p_2,q_1,q_2 \le \infty$ satisfying the relation \[ \frac{1}{r}=\frac{1}{p_1}+\frac{1}{p_2}=\frac{1}{q_1}+\frac{1}{q_2}.\]
Then, the following fractional Leibniz rule holds:
$$\|\,|D|^s(u \,v)\|_{L^r}\lesssim \|\,|D|^s u\|_{L^{p_1}}\|v\|_{L^{p_2}}+\|u\|_{L^{q_1}}\|\,|D|^s v\|_{L^{q_2}} $$
for any $u\in \dot{H}^s_{p_1} \cap L^{q_1}$ and $v\in \dot{H}^s_{q_2} \cap L^{p_2}$.
\end{md}
These results can be found in \cite{Grafakos}.
\medskip

%==================================================================
\noindent \textit{A.3. Fractional chain rule}

\begin{md} \label{Propfractionalchainrulegeneral}
Let us choose $s>0$, $p>\lceil s \rceil$
 and $1<r,r_1,r_2<\infty$ satisfying $\frac{1}{r}=\frac{p-1}{r_1}+\frac{1}{r_2}$. Let us denote by $F(u)$ one of the functions $|u|^p, \pm |u|^{p-1}u$. Then, it holds the following fractional chain rule:
$$ \|\,|D|^{s} F(u)\|_{L^r}\lesssim \|u\|_{L^{r_1}}^{p-1}\|\,|D|^{s} u\|_{L^{r_2}} $$
for any $u\in  L^{r_1} \cap \dot{H}^{s}_{r_2}$.
 \end{md}
The proof can be found in \cite{Palmierithesis}.
\medskip

%==================================================================
\noindent \textit{A.4. Fractional powers}

\begin{md} \label{PropSickelfractional}
Let $p>1$, $1< r <\infty$ and $u \in H^{s}_r$, where $s \in \big(\frac{n}{r},p\big)$.
Let us denote by $F(u)$ one of the functions $|u|^p,\, \pm |u|^{p-1}u$ with $p>1$. Then, the following estimate holds$:$
$$\Vert F(u)\Vert_{H^{s}_r}\lesssim \|u\|_{H^{s}_r}\|u\|_{L^\infty}^{p-1}.$$
\end{md}

\begin{hq} \label{Corfractionalhomogeneous}
Under the assumptions of Proposition \ref{PropSickelfractional} it holds: $\| F(u)\|_{\dot{H}^{s}_r}\lesssim \| u\|_{\dot{H}^{s}_r}\|u\|_{L^\infty}^{p-1}.$
\end{hq}
The proof can be found in \cite{DuongKainaneReissig}.
\medskip

%==================================================================
\noindent \textit{A.5. A fractional Sobolev embedding}

\bmd \label{FracSobolevEmbedding}
Let $n \ge 1$, $0< s< n$, $1< q\le r< \ity$, $\alpha< \frac{n}{q^{'}}$ where $q^{'}$ denotes conjugate number of $q$, and $\gamma > -\frac{n}{r}$, $\alpha \ge \gamma$ satisfying $\frac{1}{r}= \frac{1}{q}+ \frac{\alpha- \gamma- s}{n}$. Then, it holds:
$$ \big\||x|^\gamma |D|^{-s}u\big\|_{L^r} \lesssim \big\||x|^\alpha u\big\|_{L^q}, \text{ that is }, \big\||x|^\gamma u\big\|_{L^r} \lesssim \big\||x|^\alpha\, |D|^s u\big\|_{L^q} $$
for any $u \in \dot{H}^{s,q}_\alpha$, where $\dot{H}^{s,q}_\alpha= \{u \,:\, |D|^s u \in L^q(\R^n,\, |x|^{\alpha q})\}$ is the weighted homogeneous Sobolev space of potential type with the norm $\|u\|_{\dot{H}^{s,q}_\alpha}= \big\||x|^\alpha\, |D|^s u \big\|_{L^q}$.
\emd
The proof can be found in \cite{SteinWeiss}.

\bhq \label{CorollaryEmbedding}
Let $1< q< \ity$ and $0< s_1< \frac{n}{q}< s_2$. Then, for any function $u \in \dot{H}^{s_1,q} \cap \dot{H}^{s_2,q}$ we have
\[ \|u\|_{L^\ity} \lesssim \|u\|_{\dot{H}^{s_1,q}}+ \|u\|_{\dot{H}^{s_2,q}}. \]
\ehq
\renewcommand{\proofname}{Proof.}
\begin{proof}
By choosing $\alpha= \gamma= 0$ and $s= s_1$ in Proposition \ref{FracSobolevEmbedding} we get
$$ \|u\|_{L^r} \lesssim \big\||D|^{s_1} u\big\|_{L^q}, \text{ where } \frac{1}{r}= \frac{1}{q}- \frac{s_1}{n}. $$
Since $s_2- s_1 >\frac{n}{r}$, we can conclude
\[ \|u\|_{L^\ity} \lesssim \|u\|_{H^{s_2- s_1,r}} \lesssim \|u\|_{L^r}+ \big\||D|^{s_2- s_1} u\big\|_{L^r} \lesssim \big\||D|^{s_1} u\big\|_{L^q}+ \big\||D|^{s_2} u\big\|_{L^q}. \]
Hence, the proof of Corollary \ref{CorollaryEmbedding} is completed.
\end{proof}

%==================================================================
\noindent \textit{A.6. Modified Bessel functions}

\begin{md} \label{FourierModifiedBesselfunctions}
Let $f \in L^p(\R^n)$, $p\in [1,2]$, be a radial function. Then, the Fourier transform $F(f)$ is also a radial function and it satisfies
$$ F_n(\xi):= F(f)(\xi)= c \int_0^\ity g(r) r^{n-1} \tilde{J}_{\frac{n}{2}-1}(r|\xi|)dr,\,\,\, g(|x|):= f(x), $$
where $\tilde{J}_\mu(s):=\frac{J_\mu(s)}{s^\mu}$ is called the modified Bessel function with the Bessel function $J_\mu(s)$ and a non-negative integer $\mu$.
\end{md}

\begin{md} \label{PropertiesModifiedBesselfunctions}
The following properties hold for the modified Bessel functions:
\begin{enumerate}
\item $sd_s\tilde{J}_\mu(s)= \tilde{J}_{\mu-1}(s)-2\mu \tilde{J}_\mu(s)$,
\item $d_s\tilde{J}_\mu(s)= -s\tilde{J}_{\mu+1}(s)$,
\item $\tilde{J}_{-\frac{1}{2}}(s)= \sqrt{\frac{2}{\pi}}\cos s$ and $\tilde{J}_{\frac{1}{2}}(s)= \sqrt{\frac{2}{\pi}} \frac{\sin s}{s}$,
\item $|\tilde{J}_\mu(s)| \le Ce^{\pi|Im\mu|} \text{ if } s \le 1, $ and $\tilde{J}_\mu(s)= Cs^{-\frac{1}{2}}\cos \big( s-\frac{\mu}{2}\pi- \frac{\pi}{4} \big) +\mathcal{O}(|s|^{-\frac{3}{2}}) \text{ if } |s|\ge 1$,
\item $\tilde{J}_{\mu+1}(r|x|)= -\frac{1}{r|x|^2}\partial_r \tilde{J}_\mu(r|x|)$, $r \ne 0$, $x \ne 0$.
\end{enumerate}
\end{md}
\medskip

%==================================================================
\noindent \textit{A.7. Fa\`{a} di Bruno's formula}

\begin{md} \label{FadiBruno'sformula1}
Let $h\big(g(x)\big)= (h\circ g)(x)$ with $x\in \R$. Then, we have
$$ \frac{d^n}{dx^n}h\big(g(x)\big)= \sum \frac{n!}{m_1! 1!^{m_1} m_2! 2!^{m_2}\cdots m_n! n!^{m_n}}h^{(m_1+m_2+\cdots+m_n)}\big(g(x)\big) \prod_{j=1}^n \big( g^{(j)}(x) \big)^{m_j}, $$
where the sum is taken over all $n$- tuples of non-negative integers $(m_1,m_2,\cdots,m_n)$ satisfying the constraint of the Diophantine equation: $ 1\cdot m_1+ 2\cdot m_2+\cdots+ n\cdot m_n =n.$
\end{md}
For the proof one can see \cite{FrancescoBruno}.
\medskip

%==================================================================
\noindent \textit{A.8. Useful lemma}

\bbd \label{LemmaIntegral}
Let $\alpha, \beta \in \R.$ Then:
$$ I(t):=\int_0^t (1+t-\tau)^{-\alpha}(1+\tau)^{-\beta}d\tau \lesssim
\begin{cases}
(1+t)^{-\min\{\alpha, \beta\}} \hspace{2cm} \text{ if } \max\{\alpha, \beta\}>1,&\\
(1+t)^{-\min\{\alpha, \beta\}}\log(2+t) \hspace{0.45cm} \text{ if } \max\{\alpha, \beta\}=1,&\\
(1+t)^{1-\alpha-\beta} \hspace{2.5cm} \text{ if } \max\{\alpha, \beta\}<1.&\\
\end{cases} $$
\ebd
\renewcommand{\proofname}{Proof.}
\begin{proof}
Let us divide the interval $[0,t]$ into $[0,t/2]$ and $[t/2,t]$. It holds
\begin{align*}
&\frac{1}{2}(1+t) \le 1+t-s \le 1+t \text{ for any } s \in [0,t/2],\\
&\frac{1}{2}(1+t) \le 1+s \le 1+t \text{ for any } s \in [t/2,t].
\end{align*}
Hence, using the change of variables when needed we get
\begin{align*}
&I(t) \approx  (1+t)^{-\alpha} \int_0^{t/2}(1+\tau)^{-\beta}d\tau+ (1+t)^{-\beta} \int_{t/2}^t (1+t-\tau)^{-\alpha}d\tau\\
&\quad = (1+t)^{-\alpha} \int_0^{t/2}(1+\tau)^{-\beta}d\tau+ (1+t)^{-\beta} \int_0^{t/2}(1+\tau)^{-\alpha}d\tau \approx (1+t)^{-\min\{\alpha, \beta\}} \int_0^{t/2}(1+\tau)^{-\max\{\alpha, \beta\}}d\tau.
\end{align*}
Therefore, the proof of Lemma \ref{LemmaIntegral} is completed.
\end{proof}

%=================================================================================={References}

\end{document}